\documentclass[aap,preprint]{imsart}

\RequirePackage[OT1]{fontenc}
\RequirePackage{amsthm,amsmath}
\RequirePackage[numbers]{natbib}
\RequirePackage[colorlinks,citecolor=blue,urlcolor=blue]{hyperref}

\RequirePackage{amssymb}
\RequirePackage{graphicx}
    \graphicspath{{./pix/}} 
\RequirePackage{subfigure}

\startlocaldefs
\numberwithin{equation}{section}
\theoremstyle{plain}
\newtheorem{theorem}{theorem}[section]
\endlocaldefs

\begin{document}

\begin{frontmatter}
\title{The adaptive patched particle filter and its implementation}
\runtitle{The adaptive patched particle filter}

\begin{aug}
\author{\fnms{} \snm{Wonjung Lee}\thanksref{m1}\ead[label=e1]{leew@maths.ox.ac.uk}}
\and
\author{\fnms{} \snm{Terry Lyons}\thanksref{m1}\ead[label=e2]{tlyons@maths.ox.ac.uk}}

\runauthor{Wonjung Lee and Terry Lyons}

\affiliation{University of Oxford\thanksmark{m1}}

%

\address{
Wonjung Lee \\
Stochastic Analysis Group and OCCAM\\
Mathematical Institute \\
And\\
Oxford-Man Institute of Quantitative Finance \\
University of Oxford \\
U.K. \\
\printead{e1}\\
}

\address{
Terry Lyons\\
Stochastic Analysis Group \\
Mathematical Institute \\
And\\
Oxford-Man Institute of Quantitative Finance \\
University of Oxford \\
U.K. \\
\printead{e2}\\
}
\end{aug}

\begin{abstract}
There are numerous contexts where one wishes to describe the state of a
randomly evolving system. Effective solutions combine models that quantify
the underlying uncertainty with available observational data to form
relatively optimal estimates for the uncertainty in the system state.
Stochastic differential equations are often used to mathematically model the
underlying system. The Kusuoka-Lyons-Victoir (KLV) approach is a higher
order particle method for approximating the weak solution of a stochastic
differential equation that uses a weighted set of scenarios to approximate
the evolving probability distribution to a high order of accuracy. 
The algorithm can be performed by integrating along a number of
carefully selected bounded variation paths
and the iterated application of the KLV method has a tendency 
for the number of particles to increase.
Together with local dynamic recombination that simplifies the support of discrete measure
without harming the accuracy of the approximation,
the KLV method becomes eligible to solve the filtering problem for which 
one has to maintain an accurate description of the ever-evolving conditioned measure.
Besides the alternate application of the KLV method and recombination
for the entire family of particles,
we make use of the smooth nature of likelihood to lead some of the particles 
immediately to the next observation time and to build an algorithm that is
a form of automatic high order adaptive importance sampling.
We perform numerical simulations to evaluate the efficiency and accuracy
of the proposed approaches in the example of the linear stochastic differential equation
driven by three independent Brownian motions.
Our numerical simulations show that, even when the sequential Monte-Carlo method poorly performs, 
the KLV method and recombination can together be used to approximate higher order moments of the 
filtering solution in a moderate dimension with high accuracy and efficiency.  
\end{abstract}



\end{frontmatter}

\section{Introduction}
Filtering
is an approach for
calculating the probability distribution of an evolving system 
in the presence of
noisy observations. 
The problem
has many significant and practical applications in science and engineering,
for example
navigational and guidance systems, radar tracking, sonar ranging, 
satellite and airplane orbit determination,
the spread of hazardous plumes or pollutants,
prediction of weather and climate in atmosphere-ocean dynamics
\cite{Kalman60, 
kalman1961new,
  kushner1967approximations, 
  jazwinski1970stochastic, 
  gelb1974applied,
anderson1979optimal, 
gordon1993novel,
evensen2009data}.
If both the underlying system 
and the observation process satisfy linear equations,
the solution of the filtering problem 
can be obtained from the Kalman filter
\cite{Kalman60,
kalman1961new}.
For nonlinear filtering in finite dimension, 
there occasionally exist analytic solutions
but the results are too narrow in applicability
\cite{benevs1981exact}.
As a result,
a number of 
numerical schemes 
that use
a discrete measure,
i.e.,
collection of weighted Dirac masses,
for the approximation of 
the conditioned measure have been developed
\cite{gordon1993novel,
evensen2009data,
doucet2001sequential}.

When the underlying dynamics is a continuous process and the associated observation
is intermittent in time, 
one approach to filtering is 
to make a prediction
to quantify uncertainty
and then to update this prediction
to incorporate data
in a sequential fashion.
The prediction step
corresponds
to solving the Kolmogorov forward equation
when the system is driven by Brownian motion.
For the numerical integration of a stochastic differential equation,
the sequential Monte-Carlo method
uses sampling from a random vector
whose distribution agrees with the law of the
truncated
strong Taylor expansion of 
the solution of an Ito diffusion.
The algorithm usually gives lower order strong convergence
of the probability distribution
\cite{kloeden2011numerical}.

Instead of simulating Wiener measure as in the sequential Monte-Carlo method,
the KLV method
at the path level
replaces Brownian motion
by a weighted combination of bounded variation paths
while making sure that expectations of the 
iterated integrals 
with respect to
these two measures on Wiener space agree up to a certain degree. 
Then the particles are 
deterministically 
pushed 
forward 
along the paths
to yield a weighted discrete measure.
The KLV method 
is
of higher order with effective and transparent error bounds
obtained from
the Stratonovich-Taylor expansion of the solution of a stochastic differential equation
\cite{lyons2004cubature}.

It is intrinsic to the KLV method that the number of particles increases 
when the algorithm is iterated.
Therefore
its successive application
without an efficient suppression of the growth of the number of particles
cannot be used to filter the ever-evolving dynamics.
Given a family of test functions, one can replace the original discrete measure 
by a simpler measure with smaller support 
whose integrations against these test functions agree with those against the original measure.
Recombination achieves the reduction of particles 
in this way
using the polynomials as test functions
\cite{litterer2012high}.
One advantage of recombination is its local applicability in space.
Therefore one can divide the set of particles into 
a number of disjoint subsets and recombine 
each clustered discrete measure separately,
a process which we call the patched recombination.
The dynamic property of patched recombination,
if an efficient classification method is provided,
leads to
a competitive high order reduction algorithm
whose error bound
can be obtained from the
Taylor expansion of the test function.

One can use 
the alternate application of
the KLV method
and patched recombination 
as an algorithm for 
the prediction step
in filtering.
However the cost of this non adaptive method
would become extremely high particularly in high dimension.
In this paper we modify 
the algorithm
to reduce the computational effort.
More precisely,
we exploit the internal smoothness of the likelihood 
to allow some particles to immediately leap to the next observation time
provided the support of the resulting measure is far from the observational data.
Applying bootstrap reweighting 
to discrete measures
for the updating step,
our solution of the filtering problem
is consistent with Bayesian statistics.

It is very important to use a good example
to examine the performance of the algorithm we have developed.
We choose a forward model and observation process
for which the analytic solution of the filtering problem is known
and use this to measure the accuracy of our approximations.
Unlike filtering in practice
where the data is determined by a realization of the
random process and the observation noise,
we arbitrarily fix the observational data to study cases
ranging from  normal to exceptional and to rare event.
This setting is admittedly somewhat artificial 
however it is carefully designed in order to find the parameter regimes
where our approach outperforms Monte-Carlo methods
and eventually turns out to be extremely helpful
for a deeper understanding of the filtering problem.

The paper is organized as follows.
In section~\ref{sec:BF},
we introduce the Bayesian filter.
In section~\ref{sec:PF},
we 
describe 
a prototypical sequential Monte-Carlo algorithm and its 
variant.
Section~\ref{sec:PPF} defines the patched particle filter
and
the adaptive patched particle filter.
In section~\ref{sec:adaptive},
we discuss the implementation of 
the KLV method.
Numerical simulations are performed in section~\ref{sec:numerical}
and concluding discussions are in section~\ref{sec:discussion}.

\section{Bayesian filter}
\label{sec:BF}
Suppose that
the 
$N$-dimensional
underlying Markov process 
$X(t), t \in \mathbb{R}^+ \cup \{0\}$,
and
the
$N'$-dimensional
observation process
$Y_n, n\in \mathbb{N}\backslash \{0\}$,
associated with 
$X_n = X(n\times T)$
are given,
for some inter-observation time $T>0$.
Let $Y_{1:n'} \equiv \lbrace Y_1,\cdots,Y_{n'} \rbrace$ 
be the path of the observation process
and
$y_{1:n'} \equiv \lbrace y_1,\cdots,y_{n'} \rbrace$ 
be a generic point in the space of paths.
We define
the measure of the conditioned variable
$X_n | Y_{1:n'}$ 
by
$\pi_{n|n'}(dx_n) = \mathbb{P}(X_n\in dx_n \vert Y_{1:n'}=y_{1:n'})$.
Given $\pi_{0|0}$ which is the law of $X(0)$,
filtering aims to find $\pi_{n|n}$ for all $n \geq 1$.

This intermittent data assimilation problem can be solved by
the alternate application of the prediction,
to obtain the prior ${\pi}_{n|n-1}$ from
${\pi}_{n-1|n-1}$,
and 
the updating,
to obtain
the posterior ${\pi}_{n|n}$ from ${\pi}_{n|n-1}$.
If the transition kernel $K(dx_n |x_{n-1})$ 
and the likelihood function $g(y_n| x_n)$,
satisfying
\begin{equation*}
  \begin{split}
	\mathbb{P}(X_n \in A\vert X_{n-1}=x_{n-1}) & =\int_A K(dx_n\vert x_{n-1}), \\
	\mathbb{P}(Y_n \in B\vert X_{n}=x_{n}) &=\int_B g(y_n\vert x_{n})\,dy_n,
  \end{split}
\end{equation*}
for all
$A\in \mathcal{B}(\mathbb{R}^N)$, the Borel $\sigma$-algebra,
and
$B\in \mathcal{B}(\mathbb{R}^{N'})$,
are given,
the prediction and the updating are achieved by
\begin{align}
\pi_{n|n-1}(dx_n) 
& =\int K(dx_{n}|x_{n-1}) \pi_{n-1|n-1}(dx_{n-1}),
	\label{eq:pred}
\\
\pi_{n|n}(dx_n)
& =\frac{g\left(y_n|x_n\right) \pi_{n|n-1}(dx_n)}{\int g\left(y_n|x_n\right) \pi_{n|n-1}(dx_n)},
	\label{eq:bayes}
\end{align}
respectively.
Eq.~(\ref{eq:bayes}) is Bayes' rule and the recursive scheme 
(\ref{eq:pred}),
(\ref{eq:bayes})
is called a Bayesian filter.

\section{Particle filtering}
\label{sec:PF}
\subsection{Weak approximation}
The closed form of
$\pi_{n|n'}$ 
is in general not available.
In many cases the essential properties of a probability measure we are interested can
accurately be described 
by the expectation of 
test functions. 
If the class of test functions is specified, we can replace the
original measure with a simpler measure that integrates the test functions
correctly and hence still 
keeps the right properties of the original measure.
Therefore
efforts have been devoted to weakly
approximating $\pi_{n|n'}$
by finding an efficient way to compute
$\mathbb{E}(f(X_n) | Y_{1:n'}) = \int f(x_n) \pi_{n|n'}(dx_n)$
accurately for a sufficiently large class of 
$f:\mathbb{R}^N \to \mathbb{R}$.
We mention that the class of test functions is not given in the filtering problem.
Their choice is quite critical as it affects the notion of an optimal algorithm 
and controls the detailed description of the conditioned measure.

One of the methodologies
for the weak approximation
is to employ
\textit{particles}
whose locations and weights
characterize the approximation of the conditioned measure.
More precisely, a particle filter is a recursive algorithm that produces
\begin{equation}
  \label{eq:weighteddiscretemeasure}
  {\pi}^{\text{PF}}_{n|n'} = \sum_{i=1}^{M_{n|n'}} \lambda^i_{n|n'} \delta_{x^i_{n|n'}}
\end{equation}
approximating ${\pi}_{n|n'}$,
where $\delta_x$ denotes the Dirac mass centered at $x$.
One approximates
$(\pi_{n|n'},f)$
by
$({\pi}^{\text{PF}}_{n|n'},f) = \sum_{i=1}^{M_{n|n'}} \lambda^i_{n|n'} f(x^i_{n|n'})$
where
the notation $(\pi,f)=\int f(x) \pi(dx)$ is used.

\subsection{Sequential Monte-Carlo method}
Particle approximation is widely used in Monte-Carlo methods.
We here introduce the bootstrap filter or
sampling importance resampling (SIR) suggested in
\cite{gordon1993novel}
and the sequential importance sampling (SIS) algorithm
\cite{liu1998sequential, pitt1999filtering, doucet2000sequential}.
The number of particles do not have to be equal in each step,
but we here fix it by $M_{n|n'} = M$ for simplicity.

\subsubsection{Bootstrap filter or sampling importance resampling (SIR)}
The prediction is achieved by using
$(\pi_{n|n-1},f)=(\pi_{n-1|n-1},Kf)$
from Eq.~(\ref{eq:pred}).
Given the empirical measure
${\pi}_{n-1|n-1}^{\text{SIR}}=\frac{1}{M} \sum_{i=1}^M \delta_{x^i_{n-1|n-1}}$
approximating $\pi_{n-1|n-1}$,
one performs
independent and identically distributed (i.i.d.)
sampling
$\bar{x}^{i}_{n|n-1}$ drawn from $K(dx_n| x^{i}_{n-1|n-1})$.
Then
the discrete measure
${\pi}_{n|n-1}^{\text{SIR}} = \frac{1}{M} \sum_{i=1}^M \delta_{\bar{x}^{i}_{n|n-1}}$
is distributed 
according to 
$\pi_{n|n-1}$.

For the updating,
Eq.~(\ref{eq:bayes}) implies that
$(\pi_{n|n},f)=(\pi_{n|n-1},fg^{y_n})/(\pi_{n|n-1},g^{y_n})$
where the notation 
$g^{y_n}(\cdot) \equiv g(y_n| \cdot)$ is used.
We are led to
define the bootstrap reweighting operator
\begin{equation}
	\label{eq:bootstrapreweighting}
\text{REW}\left( \sum_{i=1}^{n} \kappa_i \delta_{x^i},g^{y_n}\right)
\equiv
\frac{\sum_{i=1}^{n} \kappa_i g^{y_n}({x}^i)\delta_{{x}^i}}{\sum_{i=1}^{n} \kappa_i g^{y_n}({x}^i)}.
\end{equation}
Then $\bar{\pi}^{\text{SIR}}_{n|n}= \text{REW}\left( \pi^{\text{SIR}}_{n|n-1},g^{y_n}\right)$
is distributed according to $\pi_{n|n}$.

In order to prevent degeneracy in the weights,
one approximates the weighted discrete measure
$\bar{\pi}^{\text{SIR}}_{n|n}$
by an equally weighted discrete measure
\cite{doucet2001sequential}.
Random resampling performs
$M$ independent samples $\{ x^i_{n|n} \}_{i=1}^M$ from $\bar{\pi}^{\text{SIR}}_{n|n}$.
This process can introduce
a large Monte-Carlo variation 
and work has been done to reduce the variance
\cite{carpenter1999improved,
crisan2002minimal}.
The resulting discrete measure
${\pi}_{n|n}^{\text{SIR}} = \frac{1}{M} \sum_{i=1}^M \delta_{{x}^{i}_{n|n}}$
is distributed according to $\pi_{n|n}$.

The SIR algorithm can be displayed by
\begin{equation}
  \label{eq:SIR}
{\pi}_{n-1|n-1}^{\text{SIR}} 
\mapsto 
{\pi}_{n|n-1}^{\text{SIR}} 
\Rightarrow
\bar{\pi}_{n|n}^{\text{SIR}} 
\to
{\pi}_{n|n}^{\text{SIR}} 
\end{equation}
where the notation $\mapsto$ is used for moving particles forward in time,
$\Rightarrow$ for reweighting and
$\to$ for random resampling.
The algorithm is very intuitive and straightforward to implement.
Further, it produces an approximation that converges toward to the true posterior
as the number of particles increases
\cite{crisan2002survey}.
However SIS might be inefficient when 
${\pi}_{n|n-1}^{\text{SIR}}$ is far from $\pi_{n|n}$ in the sense that
bootstrap reweighting generates importance weights with a high variance.
The following
SIS
algorithm modifies SIR to get around this problem.

\subsubsection{Sequential importance sampling (SIS)}
Given the unweighted measure
${\pi}_{n-1|n-1}^{\text{SIS}}=\frac{1}{M} \sum_{i=1}^M \delta_{x^i_{n-1|n-1}}$
that approximates $\pi_{n-1|n-1}$,
one performs
i.i.d. 
sampling
$\widetilde{x}^{i}_{n|n-1} \sim \widetilde{K}(dx_n| x^{i}_{n-1|n-1},y_n)$
instead of
$\bar{x}^{i}_{n|n-1} \sim K(dx_n| x^{i}_{n-1|n-1})$.
Here the new transition kernel $\widetilde{K}$ 
can depend on $y_n$
and should be chosen in a way that
the distribution of 
${\pi}_{n|n-1}^{\text{SIS}}=\frac{1}{M} \sum_{i=1}^M \delta_{\widetilde{x}^i_{n|n-1}}$
is closer to $\pi_{n|n}$ than
${\pi}_{n|n-1}^{\text{SIR}}$
in the above-mentioned sense
\cite{doucet2000sequential}.

Note that
${\pi}_{n|n-1}^{\text{SIS}}$ is not distributed according to
$\pi_{n|n-1}$.
To account for the effect of this discrepancy,
the expression
\begin{equation}
  \begin{split}
& \mathbb{P}(X_{n-1} \in dx_{n-1},X_n \in dx_{n}|Y_{1:n} = y_{1:n}) \\
& \quad = \frac{ w(x_{n-1},x_n,y_n) \widetilde{K}(dx_n|x_{n-1},y_n)\pi_{n-1|n-1}(dx_{n-1}) }
  { \int w(x_{n-1},x_n,y_n) \widetilde{K}(dx_n|x_{n-1},y_n)\pi_{n-1|n-1}(dx_{n-1}) }
  \end{split}
	\label{eq:marginal}
\end{equation}
where
\begin{equation*}
  w(x_{n-1},x_n,y_n) \propto 
  \frac{g (y_n|x_n) K(dx_n|x_{n-1})}{\widetilde{K}(dx_n|x_{n-1},y_n) }
\end{equation*}
is used.
Replacing $\widetilde{K}(dx_n|x_{n-1})\pi_{n-1|n-1}(dx_{n-1})$ in 
Eq.~(\ref{eq:marginal}) by its empirical approximation and
integrating over $x_{n-1}$,
one obtains
$\widetilde{\pi}_{n|n}^{\text{SIS}} = \sum_{i=1}^M w^i \delta_{ \widetilde{x}^i_{n|n-1} }$ 
where $w^i \propto  w(x^i_{n-1|n-1},\widetilde{x}^i_{n|n-1},y_n) $
that is distributed according to
$\pi_{n|n}$.
Random sampling from 
$\widetilde{\pi}_{n|n}^{\text{SIS}}$ 
yields the empirical measure
${\pi}_{n|n}^{\text{SIS}}$.

If $\widetilde{K}(dx_n|x_{n-1},y_n)$ and $w(x_{n-1},x_n,y_n)$ have better theoretical properties
than ${K}(dx_n|x_{n-1})$ and $g(y_n|x_n)$ such as 
better mixing properties of 
$\widetilde{K}(dx_n|x_{n-1},y_n)$ 
or flatter likelihood
$w(x_{n-1},x_n,y_n)$,
then the algorithm will outperform.
Because one needs to integrate an evolution equation of a Markov process with transition kernel 
$\widetilde{K}$ in any practical implementation, 
designing efficient particle filtering methods is equivalent to 
finding an appropriate dynamic model 
that has good theoretical properties while keeping the same filtering distributions. 
In the numerical simulations performed in 
\cite{van2010nonlinear},
the SIS algorithm
\begin{equation}
  \label{eq:SIS}
{\pi}_{n-1|n-1}^{\text{SIS}} 
\mapsto 
{\pi}_{n|n-1}^{\text{SIS}} 
\Rightarrow
\widetilde{\pi}_{n|n}^{\text{SIS}} 
\to
{\pi}_{n|n}^{\text{SIS}} 
\end{equation}
uses far fewer particles 
than SIR
to achieve a given degree of accuracy.

In the subsequent section, 
we develop two non Monte-Carlo particle filtering algorithms that 
retains the strengths and mitigates the weaknesses of the SIR and SIS methods.

\section{Introducing cubature to filtering}
\label{sec:PPF}
Suppose that
a random vector $X(t) \in \mathbb{R}^N$ 
evolves according to a Stratonovich 
stochastic differential equation (SDE)
\begin{equation}
  \label{eq:dynamics}
  dX(t)  = V_0(X(t))\,dt +\sum_{i=1}^d V_i(X(t)) \circ dW_i(t)
\end{equation}
where $\{ V_i 
\in C_b^\infty( \mathbb{R}^N,\mathbb{R}^{N})
\}_{i=0}^d$ 
is a family of smooth vector fields 
from $\mathbb{R}^N$
to $\mathbb{R}^N$
with bounded derivatives of all orders,
and $W =(W_1,\cdots, W_d)$ denote a set of Brownian motions, independent of one another.
The noisy observations $Y_n$ associated with $X_n = X(n\times T)$ 
are given by
\begin{equation}
\label{eq:observation}
Y_n  = \varphi(X_n) + \eta_n, \quad \eta_n \sim \mathcal{N}(0,R_n)  
\end{equation}
where
$\varphi \in C_b^\infty( \mathbb{R}^N,\mathbb{R}^{N'})$
and realizations of the noise $\eta_n$ are 
i.i.d. 
random vectors
in $\mathbb{R}^{N'}$.

The two main ingredients in developing 
the patched particle filter (PPF)
and the adaptive patched particle filter (APPF)
to solve the filtering problem 
(\ref{eq:dynamics}),
(\ref{eq:observation})
are
cubature on 
a finite dimensional space
and cubature on 
(infinite dimensional)
Wiener space.
Both are discrete measures 
and defined in
subsection~\ref{sec:cubonRN} and subsection~\ref{sec:cubonws},
respectively.
We describe the KLV method in subsection~\ref{sec:cubonws}
and the patched recombination in subsection~\ref{sec:recomb}.
Subsection~\ref{sec:algorithm} provides
the PPF algorithm together with error estimates
and
subsection~\ref{sec:appf} defines the APPF.

\subsection{Cubature on a finite dimensional space}
\label{sec:cubonRN}
Let $\nu$ be a (possibly unnormalized) measure on $\mathbb{R}^N$.
A discrete measure
$\widehat{\nu}^{(r)}=\sum_{j=1}^{n_r} w_j \delta_{y^j}$
is called a cubature 
(quadrature when $N=1$)
of degree $r$ 
with respect to $\nu$,
if 
$ (\nu,q)$ equals $(\widehat{\nu}^{(r)},q) = \sum_{j=1}^{n_r} w_j q(y^j)$
for all 
polynomials
$q$
whose total degree is less than or equal to $r$.
It is proved that 
a cubature $\widehat{\nu}^{(r)}$ with respect to an arbitrary positive measure $\nu$
satisfying 
$n_r \leq \binom{N+r}{r}$
exists
\cite{putinar1997note}.

Importantly,
an error bound of 
$(\nu,F)-(\widehat{\nu}^{(r)},F) \equiv (\nu-\widehat{\nu}^{(r)},F)$
for a smooth function $F:\mathbb{R}^N \to \mathbb{R}$
can be obtained from 
the Taylor expansion.
The value of 
$F$ 
at $x=(x_1,\cdots,x_N)$ is written as
\begin{equation}
  \label{eq:taylorex}
F(x) = \sum_{|\alpha| \leq r} \frac{D^\alpha F(x_0)}{\alpha!}(x-x_0)^{ \alpha}+R^r(x,x_0,F)
\end{equation}
where 
$\alpha\equiv (\alpha_1,\cdots,\alpha_N)$,
$|\alpha| \equiv \alpha_1+\cdots+\alpha_N$,
$\alpha! \equiv \alpha_1!\cdots\alpha_N!$,
$D^\alpha \equiv \partial x_1^{\alpha_1} \cdots \partial x_N^{\alpha_N}$,
$x^{ \alpha} \equiv x_1^{\alpha_1}\cdots x_N^{\alpha_N}$
and
\begin{equation}
  \label{eq:rem}
  R^r(x,x_0,F)=\sum_{|\alpha|=r+1} \frac{D^\alpha F(x^*)}{\alpha!}(x-x_0)^\alpha
\end{equation}
for some $x^* \in \mathbb{R}^N$.
If the support of $\nu$ is in a closed ball of center $x_0$ and radius $u$,
denoted by $B(x_0,u)$,
then we have
\begin{align}
  \label{eq:taylorremainder}
\lvert(\nu-\widehat{\nu}^{(r)},F) \rvert
= \lvert (\nu-\widehat{\nu}^{(r)},R^r) \rvert  
& \leq 2 (\nu,1)\parallel R^r\parallel_{L_\infty({B(x_0,u)})} \nonumber \\
& \leq \frac{ Cu^{r+1}}{(r+1)!} \sup_{|\alpha|=r+1} {\parallel D^{\alpha}F
\parallel_{L_\infty({B(x_0,u)})}}.
\end{align}
Here and after,
$C$ denotes a constant.
Eq.~(\ref{eq:taylorremainder}) reveals that
cubature on a finite dimensional space 
is an approach for numerical integration of functions on finite dimensional space
with a clear error bound.

Let $f$ be a function defined on a closed set $B \subseteq \mathbb{R}^N $.
We call $f$ by a Lipschitz function 
if there exits a constant $C'$,
$\{f^i\}_{i=0}^{\rho-1}$($f^0=f$)
and
$R_i: B \times B \to \mathbb{R}$
such that
$|f| \leq C'$, $|f^i|\leq C'$,
$f^i(x)=\sum_{l=0}^{\rho-1-i} f^{i+l}(y)(x-y)^{ l}/l! +R_i(x,y)$
and $|R_i(x,y)|\leq C' \parallel x-y \parallel^{\rho-i}$.
The smallest $C'$ for which the inequalities hold for all integer $i \in [0,\rho)$ is called 
the Lipschitz norm of $f$, denoted by 
$\parallel f\parallel_{\text{Lip}(\rho)}$.
Note that $f$ is defined locally, but can be extended to the entire space
by the Whitney theorem
\cite{stein1993singular}.
Eq.~(\ref{eq:taylorremainder}) implies
\begin{equation}
  \label{eq:Liperror}
\lvert(\nu-\widehat{\nu}^{(r)},F) \rvert
\leq \frac{Cu^{r+1} }{(r+1)!}
 {\parallel F \parallel}_{\text{Lip}(r+1)}.
\end{equation}
Note that $f$ is Lipschitz continuous if and only if 
$\parallel f \parallel_{\text{Lip}(1)}$ is finite.

\subsection{Cubature on Wiener space and the KLV method}
\label{sec:cubonws}
Consider 
the iterated integral with respect to 
$W =(W_1,\cdots, W_d)$,
\begin{equation*}
\mathcal{J}^I_{0,T} (\circ W) 
\equiv \int_{0<t_1<\cdots<t_l<T} \circ\, dW_{i_1}(t_{1}) \cdots \circ dW_{i_l}(t_{l}),
\end{equation*}
and 
the iterated integral with respect to
a continuous path of bounded variation
$\omega_T=(\omega_{T,1},\cdots,\omega_{T,d}): [0,T] \to \mathbb{R}^d$,
\begin{equation*}
\mathcal{J}^I_{0,T} (\omega_{T}) 
 \equiv \int_{0<t_1<\cdots<t_l<T} d\omega_{T,i_1}(t_{1}) \cdots d\omega_{T,i_l}(t_{l}),
\end{equation*}
where the notations
$W_0(t) = t$,
$\omega_{T,0}(t) = t$
and 
$I = (i_1,\cdots,i_l) \in \{0,\cdots,d\}^l$ 
are used.
Recall that Wiener space
$C^0_{0}\left( [0,T], \mathbb{R}^{d}\right)$ is
the set of continuous functions starting at zero.
We define
a discrete measure 
$\mathbb{Q}^{m}_T = \sum_{j=1}^{n_m} \lambda_j \delta_{\omega^{j}_{T}}$
supported on continuous paths of bounded variation 
to be a cubature on Wiener space of degree $m$
with respect to the Wiener measure
$\mathbb{P}$,
if the equation
\begin{align}
  \label{eq:momentmatching}
  \mathbb{E}_\mathbb{P} \left( \mathcal{J}^I_{0,T}(\circ W)  \right)
& = \mathbb{E}_{\mathbb{Q}^{m}_T} \left( \mathcal{J}^I_{0,T}(\circ W)  \right) \nonumber \\
& = \sum_{j=1}^{n_m} \lambda_j \mathcal{J}^I_{0,T}(\omega^{j}_{T})
\end{align}
holds for all 
$I$
satisfying
$||I|| \equiv l+\text{card}\{j,i_j=0\} \leq m$.
The existence of
$\mathbb{Q}^{m}_{T}$
with $n_m \leq \text{card} \{ I: \|I \| \leq m \}$
is proved in
\cite{lyons2004cubature}.

Similarly with the case of cubature on $\mathbb{R}^N$,
cubature on Wiener space 
can be used 
to approximate 
$\mathbb{E}_\mathbb{P}(f(X^x_{T}))$
for the random process $X^x_{t}$ in $N$ dimension 
satisfying
\begin{equation}
  \label{eq:dynamic}
dX^x_{t} = V_0(X^x_{t})\,dt + \sum_{i=1}^d V_i(X^x_{t}) \circ dW_i(t)
\end{equation}
and $X^x_{0}=x$.
The expectation of $f(X^x_T)$ against Wiener measure
can be viewed as an integral with respect to infinite dimensional Wiener space.

Let $t \mapsto X_t^{x,\omega^{j}_{\Delta}}$ for $t \in [0,\Delta]$ 
be the deterministic process satisfying
\begin{equation}
  \label{eq:nonauto}
dX_{t}^{x,\omega^{j}_\Delta} = \sum_{i=0}^d V_i(X_{t}^{x,\omega^{j}_\Delta}) \,d\omega^{j}_{\Delta,i}(t)
\end{equation}
and $X_{0}^{x,\omega^{j}_{\Delta}}=x$.
The ordinary differential equations (ODEs) of Eq.~(\ref{eq:nonauto})
are obtained from replacing the Brownian motions
$W$ 
in Eq.~(\ref{eq:dynamic})
by the bounded variation path $\omega^{j}_{\Delta}$.
The measure
$\sum_{j=1}^{n_m} \lambda_j \delta_{X_{T}^{x,\omega^{j}_{T}}}$
on $\mathbb{R}^N$ is called
the cubature approximation 
of the law of $X^x_T$ at the path level.
Note that
this discrete measure obtained from solving ODEs
is in general not a cubature with respect to the law of $X^x_{T}$.

An error estimate for 
the weak approximation of this particle method
can be derived from
the Stratonovich-Taylor expansion
of a smooth function $f$,
\begin{equation}
  \label{eq:staylor}
   f(X_T^x) =
   \sum_{||I||\leq m} V_If(x)
	\mathcal{J}^I_{0,T}(\circ W)
   +R_m(x,T,f)
\end{equation}
where the remainder $R_m(x,T,f)$ satisfies
\begin{equation}
  \label{eq:remainder}
   \sup_{x\in \mathbb{R}^N} 
    \sqrt{\mathbb{E}_{\mathbb{P}}(R_m(x,T,f)^2)} 
 \leq C \sum_{i=m+1}^{m+2} T^{i/2} \sup_{\parallel I \parallel= i}  \parallel V_If \parallel_\infty 
\end{equation}
for a constant $C$ depending on $d$ and $m$
\cite{kloeden2011numerical}.
Here
the vector field $V_i=(V_{i,1},\cdots,V_{i,N})$
is used as the differential operator
$V_i \equiv \sum_{j=1}^N V_{i,j} \partial x_j$
and 
$V_I$ denotes $V_{i_1} \cdots V_{i_k}$ 
Note that 
$\mathcal{J}^I_{0,T}(\circ W)$
works as a basis of the expansion in 
Eq.~(\ref{eq:staylor})
analogous to the monomial in 
Eq.~(\ref{eq:taylorex}).

The process $R_m(x,T,f)$ further satisfies
\begin{equation}
  \label{eq:cubremainder}
   \sup_{x\in \mathbb{R}^N} 
    {\mathbb{E}_{\mathbb{Q}^m_T}( \lvert R_m(x,T,f) \rvert )} 
 \leq C \sum_{i=m+1}^{m+2} T^{i/2} \sup_{\parallel I \parallel= i}  \parallel V_If \parallel_\infty 
\end{equation}
for a constant $C$ depending on $d$, $m$ and $\mathbb{Q}^m_1$
\cite{lyons2004cubature}.
Then the error bound of the cubature approximation at the path level
is given by
\begin{align}
	\label{eq:errorbound}
&  \sup_{x \in \mathbb{R}^N} 
  \left\lvert 
\mathbb{E}_\mathbb{P}(f(X^x_{T}))
-\sum_{j=1}^{n_m} \lambda_j f(X_{T}^{x,\omega^{j}_{T}}) 
  \right\rvert \nonumber \\
& \quad = 
\parallel  (P_T-Q^{m}_T)f \parallel_\infty
\leq C \sum_{i=m+1}^{m+2} T^{{i}/{2}}  \sup_{\parallel I \parallel= i} \parallel V_I f \parallel_\infty
\end{align}
for smooth $f$,
from
Eq.~(\ref{eq:momentmatching})
and
Eqs.~(\ref{eq:staylor}),
(\ref{eq:remainder}),
(\ref{eq:cubremainder}).
The operators 
$P_T$ and $Q^m_T$
are defined 
by
$P_{T}f(x)\equiv  \mathbb{E}_\mathbb{P}(f(X^x_{T}))$
and 
$Q^{m}_{T}f(x) \equiv \mathbb{E}_{\mathbb{Q}^{m}_T}(f(X^x_{T}))$.

The algorithm
was developed by Lyons, Victoir
\cite{lyons2004cubature},
following the work of Kusuoka
\cite{kusuoka2001approximation,
kusuoka2004approximation},
so it is referred to as
the KLV method.
Eq.~(\ref{eq:errorbound}) leads to define 
\begin{equation}
  \label{eq:klvop}
\text{KLV}^{(m)}
  \left(\sum_{i=1}^n \kappa_i \delta_{x^i}, \Delta \right)
  \equiv \sum_{i=1}^n \sum_{j=1}^{n_m} \kappa_i \lambda_j \delta_{X_{\Delta}^{x^i, \omega^{j}_{\Delta}}}
\end{equation}
that may be interpreted as a Markov operator acting on discrete measures on $\mathbb{R}^N$.

In the following, assume $T\in (0,1)$ for simplicity.
One may take a higher degree $m$ 
to achieve a given degree of accuracy
in Eq.~(\ref{eq:errorbound}).
Another method to improve the accuracy of the particle approximation 
is a successive application of the KLV operator.
Let 
$\mathcal{D} = \{ 0=t_0 < t_1 < \cdots < t_k = T \}$
be a partition of $[0,T]$ with $s_j=t_j-t_{j-1}$.
Instead of
$Q^m_{T}f(x) = ({\text{KLV}^{(m)}\left(\delta_x,T \right)},f)$,
the value of $P_Tf(x) = P_{s_1}P_{s_2} \cdots P_{s_k}f(x)$
can accurately be approximated 
by a multiple step algorithm
$Q^m_{s_1}Q^m_{s_2} \cdots Q^m_{s_k}f(x)$.

Given a discrete measure ${\mu^0}$,
we define a sequence of discrete measure by
\begin{equation}
  \begin{split}
	\label{eq:sequentialklv}
\Phi_{\mathcal{D}}^{m,0}({\mu^0}) &={\mu^0},  \\
\Phi_{\mathcal{D}}^{m,j}({\mu^0})
&=\text{KLV}^{(m)}(\Phi_{\mathcal{D}}^{m,j-1}({\mu^0}),s_j) \quad 1\leq j\leq k
  \end{split}
\end{equation}
that can be viewed as Markov chain.
The inequality
\begin{align}
  \label{eq:inequality1}
	& \left\lvert P_Tf(x) -(\Phi^{m,k}_{\mathcal{D}}(\delta_x),f) \right\rvert  \nonumber \\
	& \quad = \left\lvert \sum_{j=1}^k \left(\Phi^{m,j-1}_{\mathcal{D}}(\delta_x),P_{T-{t_{j-1}}}f \right)- 
	\left(\Phi^{m,j}_{\mathcal{D}}(\delta_x),P_{T-{t_{j}}}f \right) \right\rvert \nonumber \\
	& \quad = \left\lvert \sum_{j=1}^k \left(\Phi^{m,j-1}_{\mathcal{D}}(\delta_x),
	(P_{s_j}-Q^{m}_{s_j}) P_{T-{t_{j}}}f \right) \right\rvert \nonumber \\
	& \quad \leq \sum_{j=1}^k \parallel (P_{s_j}-Q^{m}_{s_j})P_{T-t_j}f \parallel_\infty
\end{align}
obtained from the Markovian property of the KLV operator
shows that
the total error of the repeated 
KLV application
is bounded above by the sum of 
the errors over the subintervals in the partition.
Applying Eq.~(\ref{eq:errorbound})
to estimate the upper bound of
Eq.~(\ref{eq:inequality1}),
we need $P_{T-t_j}f$ to be smooth
and it is true provided $f$ is smooth. 
In this case,
the error bound
\begin{equation*}
\sup_{x \in \mathbb{R}^N} 
\left\lvert 
P_Tf(x)
-({\Phi^{m,k}_{\mathcal{D}}(\delta_x)},f) \right\rvert 
\leq C \sum_{i=m+1}^{m+2}\sum_{j=1}^k s_j^{{i}/{2}} 
 \sup_{\parallel I \parallel = i}  \parallel V_IP_{T-t_j}f \parallel_\infty 
\end{equation*}
is obtained from
Eqs.~(\ref{eq:errorbound}), (\ref{eq:inequality1}).

The case of Lipschitz continuous $f$ is of particular interest
because 
$P_tf$ is indeed smooth 
in the direction of $\{V_i\}_{i=0}^d$
with additional conditions for these vector fields.
In the following we assume
$\{V_i\}_{i=0}^d$ satisfy the UFG
and V$0$ conditions
(see \cite{crisan2006convergence}),
then
$P_tf$ is smooth 
for a Lipschitz continuous $f$
and 
the regularity estimate
\begin{equation}
  \label{eq:ugf}
\parallel V_I P_tf \parallel_\infty 
\leq \frac{C}{t^{(\parallel I\parallel-1)/2}}
\parallel \nabla f \parallel_\infty
\end{equation}
holds
for all $t\in (0,1]$,
where $C$ is a constant independent of 
$f$
\cite{kusuoka1987applications, kusuoka2003malliavin}.
Combining Eqs.~(\ref{eq:errorbound}), (\ref{eq:inequality1})
and Eq.~(\ref{eq:ugf}),
we obtain an error estimate for the KLV method in terms of the gradient of $f$,
\begin{align}
	  \label{eq:totalerror}
& \sup_{x \in \mathbb{R}^N} \left\lvert 
P_Tf(x)
-({\Phi^{m,k}_{\mathcal{D}}(\delta_x)},f) \right\rvert \nonumber \\
& \quad
\leq C
\parallel \nabla f \parallel_\infty 
\left(
s_k^{{1}/{2}}+\sum_{i=m+1}^{m+2}\sum_{j=1}^{k-1} 
\frac{s_j^{{i}/{2}}}{(T-{t_j})^{{(i-1)}/{2}}} 
\right)
\end{align}
for the Lipschitz continuous $f$,
where $C$ is a constant independent of $k$.
Here the final term in the upper bound of
Eq.~(\ref{eq:inequality1})
is estimated by
$\parallel (P_{s_k}-Q^m_{s_k})f \parallel_\infty 
\leq \parallel P_{s_k}f-f \parallel_\infty 
+ \parallel f-Q^m_{s_k}f \parallel_\infty 
\leq Cs_k^{1/2}\parallel \nabla f \parallel_\infty$
using
the boundedness of 
$\{V_i\}_{i=0}^d$.

\begin{theorem}
Let $\mathcal{D}(\gamma) = \{ t_j\}_{j=0}^k$ be
the Kusuoka partition
\cite{kusuoka2001approximation}
given by
\begin{equation}
  \label{eq:Kpartition}
 	t_j= T\left( 1-\left( 1-\frac{j}{k}\right)^\gamma\right)
\end{equation}
then the error estimate
\begin{equation}
  \label{eq:KLVconvergece}
 \sup_{x \in \mathbb{R}^N} 
 \left\lvert 
 P_Tf(x)
-({\Phi^{m,k}_{\mathcal{D}(\gamma)}(\delta_x)},f) \right\rvert  
 \leq 
C \parallel \nabla f \parallel_\infty 
T^{1/2}
k^{-(m-1)/2}
\end{equation}
is satisfied
for a Lipschitz continuous $f$
when $\gamma>m-1$.
\end{theorem}

Eq.~(\ref{eq:KLVconvergece}) is
obtained from substituting 
the non-equidistant time discretization
$\mathcal{D}(\gamma)$
into
Eq.~(\ref{eq:totalerror}).
Using this particular choice of partition ensures that the bound of the KLV error is of high order
in the number of iterations $k$.

Before concluding this subsection,
we here mention that
$u(x,t) \equiv \mathbb{E}_\mathbb{P}(f(X^x_{T-t}))$
satisfies the partial differential equation (PDE)
\begin{equation}
  \begin{split}
  \label{eq:parabolicpde}
  \frac{\partial}{\partial t}u(x,t) & = - \left(V_0 + \frac{1}{2}\sum_{i=1}^d V_i^2 \right)u(x,t), \\
  u(x,T) & =f(x).
  \end{split}
\end{equation}
where $\{V_i\}_{i=0}^d$ are used as differential operators
\cite{watanabe1981stochastic}.
Therefore $P_Tf(x)$, the heat kernel applied to $f$, is equal 
to the solution $u(x,0)$ of
Eq.~(\ref{eq:parabolicpde}).
Due to this concrete relationship between parabolic PDEs and SDEs,
one can use any well-known algorithm for the solution of 
Eq.~(\ref{eq:parabolicpde})
in the prediction step of the filtering problem
determined by Eq.~(\ref{eq:dynamics}) and Eq.~(\ref{eq:observation}).
However
it is very important to note the difference between these two problems.
One needs to weakly approximate 
the law of $X(T)$,
when $X(0)$ is given by $\delta_x$,
that accurately integrate the test function $f$
for the PDE problem
while the filtering problem requires one to approximate the posterior measure of
$X_{n}|Y_{1:n}$
for all $n \geq 1$,
for which the test function 
(and the law of $X(0)$ as well in practice)
is not specified.

\subsection{Local dynamic recombination}
\label{sec:recomb}
A successive application of KLV operator gives rise to geometric growth of the number of
particles in view of Eq.~(\ref{eq:klvop}).
Except some cases of PDE problem
in which the KLV method can produce an accurate approximation
with small number of iterations,
this geometric growth of particle number
prohibits an application of the KLV method
especially for the filtering problem
where to maintain an accurate description of the ever-evolving measure
with reasonable computational cost is an important issue.
It is therefore necessary to add a simplification of the discrete measure
to the KLV method as a way to control the number of particles in the approximation 
between the successive iterations.

Though this simplification problem can be 
solved by random resampling used in the bootstrap filter,
we here apply recombination to efficiently reduce the support of discrete measure
(see \cite{litterer2012high} for the detailed algorithm).
The method produces a measure with reduced support which preserves
the expectations of the polynomials.
In this case,
the reduced measure is a cubature on $\mathbb{R}^N$ with respect to the original measure.
Because the measure from the KLV method is used to integrate $P_tf$ 
which is smooth for Lipschitz continuous $f$,
one can use the Taylor expansion for the error estimate.

Instead of using a cubature of higher degree 
to recombine the entire family of particles all at once,
we follow the work in 
\cite{litterer2012high}
to improve the performance
by dividing a given discrete measure into locally supported unnormalized measures
and replacing each separated measure by the cubature of lower degree.
We believe that 
this local dynamic recombination is a competitive algorithm with general applications
mainly because each reduction can be performed in a parallel manner to save computational time
and the error bound from the Taylor approximation is of higher order.

Let $U=(U_i)_{i=1}^R$ be a collection of balls of radius $u$
that covers the support of discrete measure $\mu$ on $\mathbb{R}^N$,
then one can find unnormalized measures $(\mu_i)_{i=1}^R$ such that $\mu=\bigsqcup_{i=1}^R \mu_i$
($\mu_i$ and $\mu_j$ have disjoint support for $i\neq j$)
and $\text{supp}(\mu_i)\subseteq U_i \cap \text{supp}(\mu)$.
In this case, 
we define the patched recombination operator by
\begin{equation}
  \label{eq:recomb}
\text{REC}^{(u,r)}
\left( \mu \right)
\equiv
\bigsqcup_{i=1}^R 
\widehat{\mu}^{(r)}_i
\end{equation}
where
$\widehat{\mu}^{(r)}_i$ denotes a cubature of degree $r$ with respect to $\mu_i$.

Given a discrete measure ${\mu^0}$,
we define a sequence of discrete measure by
\begin{equation}
  \begin{split}
  \label{eq:sequenceofdiscretemeasure}
\Phi^{m,0}_{\mathcal{D},(u,r)}({\mu^0}) & ={\mu^0}, \\
\widehat{\Phi}_{\mathcal{D},(u,r)}^{m,j-1}({\mu^0}) & = 
\text{REC}^{(u_{j-1},r_{j-1})}\left( {\Phi}_{\mathcal{D},(u,r)}^{m,j-1}({\mu^0}) \right), \\
\Phi_{\mathcal{D},(u,r)}^{m,j}({\mu^0}) &= \text{KLV}^{(m)}
\left(\widehat{\Phi}_{\mathcal{D},(u,r)}^{m,j-1}({\mu^0}),s_j \right),
  \end{split}
\end{equation}
for $1 \leq  j \leq k$.
An application of Eq.~(\ref{eq:sequenceofdiscretemeasure})
with initial condition $\delta_x$
yields a weak approximation for the law of $X^x_T$.
One obtains the estimate
\begin{align}
	  \label{eq:errorineq}
& 
\left\lvert P_Tf(x) -(\Phi^{m,k}_{\mathcal{D},(u,r)}(\delta_x),f) \right\rvert  \nonumber \\
& = \Bigg\lvert \sum_{j=1}^k
\left(\widehat{\Phi}^{m,j-1}_{\mathcal{D},(u,r)}(\delta_x),P_{T-{t_{j-1}}}f\right)- 
\left({\Phi}^{m,j}_{\mathcal{D},(u,r)}(\delta_x),P_{T-{t_{j}}}f \right) \nonumber \\
& \quad\,\, + \left({\Phi}^{m,j-1}_{\mathcal{D},(u,r)}(\delta_x),P_{T-{t_{j-1}}}f \right)- 
\left(\widehat{\Phi}^{m,j-1}_{\mathcal{D},(u,r)}(\delta_x),P_{T-{t_{j-1}}}f \right) \Bigg\rvert
\nonumber \\
& = \Bigg\lvert \sum_{j=1}^k
\left(\widehat{\Phi}^{m,j-1}_{\mathcal{D},(u,r)}(\delta_x),
(P_{s_j}-Q^m_{s_j}) P_{T-{t_{j}}}f \right)  \nonumber\\
& \quad\,\, + \left({\Phi}_{\mathcal{D},(u,r)}^{m,j-1}(\delta_x)
-\widehat{\Phi}_{\mathcal{D},(u,r)}^{m,j-1}(\delta_x),P_{T-t_{j-1}}f \right)
  \Bigg\rvert \nonumber\\
& \leq \sum_{j=1}^k \parallel (P_{s_j}-Q^m_{s_j})P_{T-t_j}f \parallel_\infty\nonumber\\
 & \quad \,\, +\sum_{j=0}^{k-1} 
 \left\lvert  
\left({\Phi}_{\mathcal{D},(u,r)}^{m,j}(\delta_x)
-\widehat{\Phi}_{\mathcal{D},(u,r)}^{m,j}(\delta_x),P_{T-t_{j}}f\right)
 \right\rvert 
\end{align}
where the first sum of the upper bound is
due to the KLV approximation.
The second sum is the error introduced by the recombination

Suppose that $f$ is Lipschitz continuous.
The smoothness of $P_tf$ leads to
\begin{align}
	\label{eq:errorbound3}
 \sup_{x \in \mathbb{R}^N} 
& \left\lvert 
\left({\Phi^{m,j}_{\mathcal{D},(u,r)}(\delta_x)}
-\widehat{\Phi}^{m,j}_{\mathcal{D},(u,r)}(\delta_x),P_{T-t_j}f \right)
\right\rvert  \nonumber \\
& \leq C u_j^{r_j+1} \sup_{|\alpha|=r_j+1} \parallel D^{\alpha}P_{T-t_j}f \parallel_\infty
\end{align}
for $0 \leq j \leq k-1$,
where
Eq.~(\ref{eq:taylorremainder}) 
and the triangle inequality are used.
In the following we assume
$\{V_i\}_{i=0}^d$ satisfy the UH condition 
(see \cite{litterer2012high, kusuoka2003malliavin}),
then
there exists a positive integer $p$ 
such that
\begin{equation}
  \label{eq:uh}
  \sup_{|\alpha|=r+1}
  \parallel
  D^\alpha P_tf
  \parallel_\infty
  \leq C t^{-rp/2}
  \parallel
  \nabla f
  \parallel_\infty
\end{equation}
for all $t\in (0,1]$.
Since the UH condition implies the UFG and V$0$ conditions
\cite{crisan2006convergence},
one obtains
\begin{align}
	\label{eq:errorbound4}
\sup_{x \in \mathbb{R}^N} 
& \left\lvert 
P_Tf(x)
-\left( {\Phi^{m,k}_{\mathcal{D},(u,r)}(\delta_x)},f \right) \right\rvert \nonumber \\
& 
\leq 
\Bigg( C_1
\Bigg(
s_k^{{1}/{2}}+\sum_{i=m+1}^{m+2}\sum_{j=1}^{k-1} 
\frac{s_j^{{i}/{2}}}{(T-{t_j})^{{(i-1)}/{2}}} 
\Bigg)  \nonumber \\
& \qquad \quad + C_2
\sum_{j=1}^{k-1} 
\frac{u_j^{r_j+1}}{(T-t_j)^{{r_jp}/{2}}}  \Bigg)
\parallel \nabla f\parallel_\infty 
\end{align}
from Eqs.~(\ref{eq:totalerror}), (\ref{eq:errorbound3}).
Here $C_1$ and $C_2$ are constants.

The recombination error can be controlled by
the radius of the ball $u_j$
and the cubature on $\mathbb{R}^N$ degree $r_j$.
By choosing a suitable pair $(u_j,r_j)$, 
one can make the order of the recombination error bound
not dominant over the order of the error bound in the KLV method.

\begin{theorem}
In the case of
$(u_j,r_j)=(s_j^{p/2-a},\lceil m/p \rceil)$
where $a = (p-1)/(2( \lceil m/p \rceil+1 ))$
($\lceil x \rceil$ denotes the smallest integer greater than or equal to $x$)
or
$(u_j,r_j)=( ( s_j^{m+1}/ (T-t_j)^{m-rp} )^{1/2(r+1)}, m)$,
the error estimate
\begin{equation}
\label{eq:klvconv}
 \sup_{x \in \mathbb{R}^N} 
\left\lvert  
P_Tf(x)
-\left( {\Phi^{m,k}_{\mathcal{D}(\gamma),(u,r)}(\delta_x)},f \right) \right\rvert 
\leq C \parallel \nabla f \parallel_\infty T^{1/2} k^{-(m-1)/2}  
\end{equation}
is satisfied
for a Lipschitz continuous $f$
when $\gamma > m-1$.
\end{theorem}

Eq.~(\ref{eq:klvconv}) is obtained from substituting 
the partition defined in Eq.~(\ref{eq:Kpartition})
into Eq.~(\ref{eq:errorbound4}).
It ensures that the recombination can be used without harming the accuracy of the KLV approximation.

\subsection{Patched particle filter}
\label{sec:algorithm}
Let $\pi_{n|n'}$ be the law of the conditioned variable $X_n|Y_{1:n'}$
determined by
Eqs.~(\ref{eq:dynamics}), (\ref{eq:observation}).
Let $ {\pi}^{\text{PPF}}_{0|0}$ be a discrete measure 
distributed according to
the law of $X_0$.
We define the 
\textit{patched particle filter (PPF) at the path level}
by the recursive algorithm
\begin{equation}
  \begin{split}
	\label{eq:patchedparticlefitler}
{\pi}^{\text{PPF}}_{n|n-1} & =
\Phi^{m,k}_{\mathcal{D},(u,r)}({\pi}^{\text{PPF}}_{n-1|n-1}), \\
{\pi}^{\text{PPF}}_{n|n} & = \text{REW}\left( {\pi}^{\text{PPF}}_{n|n-1}, g^{y_n}
\right),
  \end{split}
\end{equation}
for $n \geq 1$.
Recall that the bootstrap reweighting operator 
REW
is defined in 
Eq.~(\ref{eq:bootstrapreweighting}).
The PPF 
does not require random resampling
and therefore no Monte-Carlo variation is introduced.
The algorithm 
can be stated as the following.

\begin{enumerate}
  \item One breaks the measure into patches and performs recombination for each one.
  \item One moves given discrete measure forward in time using the KLV method.
  \item One performs data assimilation via bootstrap reweighting at 
	every inter-observation time 
	which might differ from the time step 
	for the numerical integration.
  \item One again applies the patched recombination.
\end{enumerate}

Using ${\pi}^{\text{PPF}}_{n-1|n-1}$ in place of $\delta_x$ in 
Eq.~(\ref{eq:errorineq}),
an error bound of the prior approximation
of the PPF
is given by
\begin{align}
	  \label{eq:priorerrorineq}
\left\lvert (\pi_{n|n-1} - {\pi}^{\text{PPF}}_{n|n-1},f ) \right\rvert  
& \leq \left\lvert  (\pi_{n-1|n-1},P_Tf )-
({\pi}^{\text{PPF}}_{n-1|n-1},P_Tf )\right\rvert \nonumber\\
& \quad + \left\lvert({\pi}^{\text{PPF}}_{n-1|n-1},P_Tf)
-(\Phi^{m,k}_{\mathcal{D},(u,r)}({\pi}^{\text{PPF}}_{n-1|n-1} ),f ) \right\rvert \nonumber\\
& \leq \left\lvert  (\pi_{n-1|n-1}- {\pi}^{\text{PPF}}_{n-1|n-1},P_Tf )\right\rvert \nonumber\\
& \quad + \sum_{j=1}^k \parallel (P_{s_j}-Q^m_{s_j}) P_{T-t_j}f \parallel_\infty\nonumber\\
 & +\sum_{j=0}^{k-1} 
 \left\lvert  
\left({\Phi}_{\mathcal{D},(u,r)}^{m,j}({\pi}^{\text{PPF}}_{n-1|n-1})
-{\widehat{\Phi}}_{\mathcal{D},(u,r)}^{m,j}({\pi}^{\text{PPF}}_{n-1|n-1}),P_{T-t_{j}}f\right)
 \right\rvert.
\end{align}
One can use the same argument with the case of PDE problem
to obtain a higher order estimate of the PPF.
An error bound of the posterior approximation 
\begin{align}
	\label{eq:posteriorestimate}
	\left\lvert (\pi_{n|n}-{\pi}^{\text{PPF}}_{n|n},f) \right\rvert 
	& =
	\Bigg\vert 
	\frac{(\pi_{n|n-1},fg^{y_n})}{({\pi}_{n|n-1},g^{y_n})}
	-
	\frac{({\pi}^{\text{PPF}}_{n|n-1},fg^{y_n})}{({\pi}_{n|n-1},g^{y_n})} \nonumber \\
	& \qquad  +
	\frac{({\pi}^{\text{PPF}}_{n|n-1},fg^{y_n})}{({\pi}_{n|n-1},g^{y_n})}
	-
	\frac{({\pi}^{\text{PPF}}_{n|n-1},fg^{y_n})}{({\pi}^{\text{PPF}}_{n|n-1},g^{y_n})}
	\Bigg\vert\nonumber \\
	&  \leq
	\frac{1}{(\pi_{n|n-1},g^{y_n})}
	\left\lvert (\pi_{n|n-1}-{\pi}^{\text{PPF}}_{n|n-1},fg^{y_n})\right\rvert  \nonumber \\
	& \qquad + \frac{\parallel f \parallel_\infty}{(\pi_{n|n-1},g^{y_n})}
	\left\lvert (\pi_{n|n-1}-{\pi}^{\text{PPF}}_{n|n-1},g^{y_n})\right\rvert  
\end{align}
is given in terms of an error estimate of the prior approximation.
Eq.~(\ref{eq:priorerrorineq})
and
Eq.~(\ref{eq:posteriorestimate})
implies
the higher order weak convergence of the PPF.

The implementation of the PPF 
requires to specify
the time partition and 
the way of dividing the support of measure into patches.
Before presenting the ones used in our numerical simulations,
we build a modification of the PPF.

\subsection{Adaptive patched particle filter}
\label{sec:appf}
If we know nothing about the smoothness of the test function 
or any objective function to integrate, 
then the accuracy requirement leaves no choice other than 
to let the cloud of particles increase at an appropriate rate 
to the observation time and the KLV numerical approach indeed allows that to happen. 
For truly irregular test functions, 
the computation would be necessarily expensive. 
This is no surprise since accurate integration would require exploration of the irregularities.

However in many settings the test function is actually
piecewise smooth and the less regular set is of significantly lower
dimension than the main part of the smoothness. 
In this case we can apply the
same analysis of smoothing but now we observe that the test function in
front of a particular point $(x,T)$ is actually far smoother that would be the case in general. 
If it is smooth enough, then we can often use our high order method
to go straight to the next observation time from some considerable distance back instead
of the step predicted in the worst case which we would otherwise have used
to terminate the algorithm.

We build this insight into the practical algorithm. 
At each application of the KLV operator,
the algorithm evaluates the test function directly using a 
one step prediction to the next observation time
and compares this with the evaluation using 
a two (or three to break certain pathological symmetries) step prediction.
If two values agree within the error tolerance, 
then the particles immediately leap to the 
next observation time.
Otherwise the prediction will follow the original partition.

In terms of accuracy,
the approach is pragmatically rather
successful because the chances of two or three steps producing consistent
answers by chance is essentially negligible. 
Furthermore,
the adaptive switch 
for which the KLV is employed
to move the prediction measure forward but
move a part of it straight to the observation time whenever the relevant part of the
test function in front of the point is smooth enough
has a very significant effect of 
pruning the computation and speeding up the algorithm
due to the reduction of particles to be recombined at each iteration.

This adaptive KLV method is an automated form of high order importance sampling and
cannot be applied without a test function. 
Differently from the PDE problem,
the test function is not specified in the filtering problem.
In practice we take the likelihood as test function to lead an adaptation.

Recall
$\mathcal{D} = \{ 0=t_0 < t_1 < \cdots < t_k = T \}$
is a partition of $[0,T]$ with $s_j=t_j-t_{j-1}$.
We use the likelihood 
$g^{y_n}$ 
to define the splitting operator
acting on
a discrete measure ${\mu}^{j-1} = \sum_{i=1}^n \kappa_i \delta_{x^i}$
at time $t_{j-1}$.
Let
${\mu}^{j-1}_{i,21}=\text{KLV}\left(\delta_{x^i},t_j-t_{j-1}\right)$,
${\mu}^{j-1}_{i,22}=\text{KLV} \left( {\mu}_{i,21}, t_k-t_j\right)$
and
${\mu}^{j-1}_{i,1}= \text{KLV}\left(\delta_{x^i},t_k-t_{j-1}\right)$.
Let $I_\tau$ be the collection of index $i$ satisfying
$\vert  ({\mu}^{j-1}_{i,1}- {\mu}^{j-1}_{i,22},  g^{y_n}) \vert<\tau$.
Then the discrete measure ${\mu}^{j-1}$ is the union of two discrete measure
$ \mu^{j-1} = \mu^{j-1, < \tau} \sqcup \mu^{j-1, \geq \tau}$
where
$\mu^{j-1, < \tau} =  \sum_{i \in I_\tau} \kappa_i \delta_{x^i}$.
For simplicity,
$\mu^{k-1, \geq \tau}$ is defined to be the null set.
The process defines the splitting operator
\begin{equation}
  \label{eq:sorting}
\text{SPL}^{(\tau)}
\left( \mu^{j-1}, g^{y_n} \right)
\equiv
\mu^{j-1, < \tau}
\end{equation}
for $1 \leq  j \leq k$.

Define a sequence of discrete measures as follows
\begin{equation}
  \begin{split}
  \label{eq:adaptivesequenceofdiscretemeasure}
\Phi^{m,0}_{\mathcal{D},(u,r),\tau}({\mu}^0) & ={\mu}^0, \\
\widehat{\Phi}_{\mathcal{D},(u,r),\tau}^{m,j-1}({\mu}^0) & = 
\text{REC}^{(u_{j-1},r_{j-1})}\left(
{\Phi}_{\mathcal{D},(u,r),\tau}^{m,j-1}({\mu}^0) \right), \\
\widehat{\Phi}^{m,j-1,< \tau}_{\mathcal{D},(u,r),\tau}({\mu}^0) & = 
\text{SPL}^{(\tau)}\left(
\widehat{\Phi}^{m,j-1}_{\mathcal{D},(u,r),\tau}({\mu}^0), g^{y_n}\right), \\
\Phi_{\mathcal{D},(u,r),\tau}^{m,j}({\mu}^0) &= \text{KLV}^{(m)}
\left(\widehat{\Phi}_{\mathcal{D},(u,r),\tau}^{m,j-1 , < \tau }({\mu}^0),s_j \right).
  \end{split}
\end{equation}
for $1 \leq  j \leq k$.
Let 
$\widehat{\Phi}^{m,j-1}_{\mathcal{D},(u,r),\tau}({\mu}^0)
=
\widehat{\Phi}^{m,j-1,< \tau}_{\mathcal{D},(u,r),\tau}({\mu}^0)
\sqcup
\widehat{\Phi}^{m,j-1, \geq \tau}_{\mathcal{D},(u,r),\tau}({\mu}^0)$
and
\begin{equation}
\Psi_{\mathcal{D},(u,r),\tau}^{m,j-1,k}({\mu}^0) 
 = \text{KLV}^{(m)} \Big( 
\text{KLV}^{(m)} 
\left( \widehat{\Phi}_{\mathcal{D},(u,r),\tau}^{m,j-1, \geq \tau}({\mu}^0),
t_j-t_{j-1} \right),T-t_j \Big).
\end{equation}
for $1 \leq  j \leq k-1$.
The \textit{adaptive patched particle filter (APPF) at the path level} is defined by
\begin{equation}
  \begin{split}
	\label{eq:adaptivepatchedparticlefitler}
{\pi}^{\text{APPF}}_{n|n-1} 
& = \left(\bigsqcup_{j=1}^{k-1} \Psi_{\mathcal{D},(u,r),\tau}^{m,j-1,k}
({{\pi}^{\text{APPF}}_{n-1|n-1}})\right), \\
& \qquad \sqcup \Phi_{\mathcal{D},(u,r),\tau}^{m,k}({{\pi}^{\text{APPF}}_{n-1|n-1}}) \\
{\pi}^{\text{APPF}}_{n|n} 
& = \text{REW}\left( {\pi}^{\text{APPF}}_{n|n-1}, g^{y_n} \right),
  \end{split}
\end{equation}
for $n \geq 1$.
The algorithm 
can be stated as the following.
\begin{enumerate}
  \item One breaks the measure into patches and performs recombination for each one.
  \item One splits given discrete measure to lead some of the particles to the next observation time
	and the rest particles to the next iteration time using the KLV method.
  \item One performs data assimilation via bootstrap reweighting at 
	every inter-observation time 
	which might differ from the time step 
	for the numerical integration.
  \item One again applies the patched recombination.
\end{enumerate}

We here mention that the likelihood is a natural choice
in view of Eq.~(\ref{eq:posteriorestimate})
for the filtering problem 
in which
the posterior is of primary interest.
We also mention that
one can apply $g^{y_n}$ and $fg^{y_n}$
simultaneously
as the test function
for the ADA operator
in Eq.~(\ref{eq:adaptivesequenceofdiscretemeasure})
if one would like to obtain a posterior approximation that
accurately integrates $f$.

It would be interesting to compare
the PPF (\ref{eq:patchedparticlefitler})
and 
SIR (\ref{eq:SIR}),
and to compare
the APPF (\ref{eq:adaptivepatchedparticlefitler})
and 
SIS (\ref{eq:SIS}).
Both PPF and SIR achieve prior approximations without using the observational data,
and subsequently achieve posterior approximations via bootstrap reweighting.
In that the observation is used to 
move particles forward in time,
the APPF is very much like 
the SIS algorithm.
However,
the APPF does not introduce new dynamics
and approximates the prior measure of the given forward model.
The way of modifying original algorithm
is different but the philosophy is the same - 
making use of the observational information
to lead the particles for an efficient approximation of the posterior.

\section{Practical implementation}
\label{sec:adaptive}
In this section we discuss several issues 
related to the implementation of the 
PPF and APPF.
Some further considerations of cubature on Wiener space
are gathered in
subsection~\ref{sec:flows}.
We make use of the test function to define an adaptive partition
(subsection~\ref{sec:adapart})
and adaptive recombination
(subsection~\ref{sec:adarecomb}).

\subsection{Cubature on Wiener space continued}
\label{sec:flows}
We here 
study the construction of cubature formula
$\mathbb{Q}^{m}_T$.
Cubature on Wiener space in terms of Lie polynomial
is defined
and 
used to develop an approximation 
based on the autonomous ODEs at flow level.

Let $\{e_i\}_{i=0}^d$
be the standard basis of $\mathbb{R} \oplus \mathbb{R}^d$.
Let $\mathcal{T}$ denote the associative and non-commutative tensor algebra 
of polynomials
generated by $\{e_i\}_{i=0}^d$.
The exponential and logarithm on $\mathcal{T}$ are defined by
\begin{equation}
  \begin{split}
	\label{eq:explog}
	\text{exp}(a) & \equiv \sum_{i = 0}^\infty \frac{{a}^{\otimes i}}{i !}, \\
  \text{log}(a) & \equiv \log(a_0)+\sum_{i = 1}^\infty
  \frac{(-1)^{i-1}}{i} \left(\frac{a}{a_0} -1 \right)^{\otimes i},
  \end{split}
\end{equation}
where $a=\sum_{I} a_I e_I$
and
$ e_I = e_{i_1}\otimes \cdots \otimes e_{i_l}$
for a multi-index
$I = (i_1,\cdots,i_l) \in \{0,\cdots,d\}^l$.
Here $\otimes$ denotes the tensor product.
We define the operators
$\text{exp}^{(m)}(\cdot)$ and $\text{log}^{(m)}(\cdot)$
by the truncation
of Eq.~(\ref{eq:explog})
to the case $\parallel I \parallel \leq m$.

We define the signature of a continuous path of bounded variation
$\omega_T:[0,T]\to \mathbb{R}^d$ by
\begin{equation*}
  \begin{split}
\mathcal{S}_{0,T}(\omega_T)
& \equiv \sum_{l=0}^\infty \int_{0<t_1<\cdots <t_l<T} d\omega_{T}(t_1)\otimes\cdots\otimes
d\omega_{T}(t_l)\\
& =  \sum_{I}   
\mathcal{J}^I_{0,T} (\omega_{T})\,e_I  
  \end{split}
\end{equation*}
and similarly the signature of a Brownian motion $W$ by
\begin{equation*}
\mathcal{S}_{0,T}(\circ W) \equiv \sum_{I }   \mathcal{J}^I_{0,T} (\circ W) \,e_I.
\end{equation*}
In view of Eq.~(\ref{eq:momentmatching}),
the definition of cubature on Wiener space of degree $m$ can be rephrased by
\begin{equation}
\label{eq:signaturemomentmatching}
\mathbb{E}_{\mathbb{P}}\left( \mathcal{S}^{(m)}_{0,T}(\circ W)\right)
=
\mathbb{E}_{\mathbb{Q}_T^m}\left( \mathcal{S}^{(m)}_{0,T}(\circ W)\right)
\end{equation}
where
$\mathcal{S}^{(m)}_{0,T}(\cdot)$ 
is the truncation of 
$\mathcal{S}_{0,T}(\cdot)$ 
to the case $\parallel I \parallel \leq m$.

Define $\mathcal{L}$ to be the space of Lie polynomials, i.e.,
linear combinations of finite sequences of Lie brackets
$[e_i,e_j]=e_i\otimes e_j - e_j \otimes e_i$.
Chen's theorem implies that 
\begin{equation}
  \label{eq:liepoly}
\mathcal{L}_T^{j} \equiv \text{log}^{(m)}( \mathcal{S}_{0,T}(\omega^{j}_{T}) )
\end{equation}
is an element of $\mathcal{L}$, i.e., a Lie polynomial.
Then the measure 
$\widetilde{\mathbb{Q}}^{m}_{T} = \sum_{j=1}^{n_m} \lambda_j \delta_{\mathcal{L}^{j}_T}$
supported on Lie polynomials
satisfies
\begin{align}
\label{eq:liepolysignaturemomentmatching}
	  \mathbb{E}_{\mathbb{P}}\left( S^{(m)}_{0,T}(\circ W)\right)
	  &= \mathbb{E}_{\widetilde{\mathbb{Q}}^m_T}\left( \text{exp}^{(m)}(\mathcal{L})  \right)
	  \nonumber \\
	  &= \sum_{j=1}^{n_m} \lambda_j \text{exp}^{(m)}(\mathcal{L}_T^{j}).
\end{align}
Conversely, for any Lie polynomials $\mathcal{L}_T^j$, 
there exists continuous bounded variation paths $\omega_T^j$
whose truncated logarithmic signature is $\mathcal{L}_T^j$.
Moreover if 
$\widetilde{\mathbb{Q}}^{m}_{T}$ satisfies
Eq.~(\ref{eq:liepolysignaturemomentmatching}),
then
${\mathbb{Q}}^{m}_{T}$ satisfies
Eq.~(\ref{eq:signaturemomentmatching}).
Therefore 
Eq.~(\ref{eq:signaturemomentmatching})
and
Eq.~(\ref{eq:liepolysignaturemomentmatching})
are equivalent.
The discrete measure
$\widetilde{\mathbb{Q}}^{m}_{T}$ is also defined as 
cubature on Wiener space.

The expectation of the truncated Brownian signature is
\begin{equation}
  \label{eq:expsig}
\mathbb{E}_{\mathbb{P}}\left( \mathcal{S}^{(m)}_{0,1}(\circ W) \right) = 
\text{exp}^{(m)}\left( e_0 + 
\frac{1}{2}
\sum_{i=1}^d e_i\otimes e_i \right)
\end{equation}
which is proved
in \cite{lyons2004cubature}.
It is immediate from Eq.~(\ref{eq:expsig}) that
cubature formulae on Wiener space for $m=2n-1$ and $m=2 n$ are equal to one another.
A formula
$\{\lambda_j,\mathcal{L}^j_1\}_{j=1}^{n_m}$
satisfying Eq.~(\ref{eq:liepolysignaturemomentmatching})
is found when $m=3$ and $m=5$ for any $d$
\cite{lyons2004cubature}.
In some cases of $m \geq 7$,
cubature formula of Lie polynomial is available when $d=1, 2$
(See \cite{gyurko2011efficient}).

From 
this $\widetilde{\mathbb{Q}}^{m}_{1}$
and 
Eq.~(\ref{eq:liepoly}),
one can construct 
${\mathbb{Q}}^{m}_{1}$
(See 
\cite{lyons2004cubature,
gyurko2011efficient}).
It follows from the scaling property of the Brownian motion that
$\omega^{j}_{T,0}(t) = \omega^{j}_{1,0}(t)$ and
$\omega^{j}_{T,i}(t) = \sqrt{T}\omega^{j}_{1,i}(t/T)$ for $1 \leq i \leq d$.
The paths define
a cubature formula $\widetilde{\mathbb{Q}}^{m}_{T}$.
Using
$\mathcal{J}^I_{0,T} (\circ W) \triangleq T^{\parallel I\parallel/2}\mathcal{J}^I_{0,1} (\circ W)$
and Eq.~(\ref{eq:liepoly}),
the scaling of the Lie polynomial 
is 
$\mathcal{L}^{j}_{T} = \langle T, \mathcal{L}^{j}_{1} \rangle$
where
$\left\langle t, \sum_{I} a_I e_I \right\rangle \equiv \sum_{I} t^{\parallel I \parallel/2}a_I e_I$.
The Lie polynomials define
a cubature formula $\mathbb{Q}^m_T$.

We next study 
the approximation
based on the flows of autonomous ODEs.
It is in fact
corresponds to a version of Kusuoka's algorithm
\cite{kusuoka2001approximation}.
Let $\Gamma$ denote the algebra homomorphism generated by $\Gamma(e_i)= V_i$ for
$i=0,\cdots,d$. 
For a vector field 
$V \in C_b^\infty( \mathbb{R}^N,\mathbb{R}^{N})$,
we define 
the flow $\text{Exp}\left(tV\right)(x) \equiv \xi_t^x$
to be 
the solution of the 
ODE
$  {d\xi_t^x} = V(\xi_t^x)\,dt$
with $\xi_0^x=x$.
By interchanging the algebra homomorphism $\Gamma$ with the exponentiation (so far taken in the
tensor algebra) we arrive at an approximation operator in which the exponentiation is understood as
taking the flow of autonomous ODEs.
More precisely,
one has
\begin{equation*}
  \begin{split}
\mathbb{E}_{\mathbb{P}}\left( \Gamma \left(S^{(m)}_{0,T}(\circ W) \right)\right) f(x)
& = \sum_{j=1}^{n_m} \lambda_j \Gamma \left(\exp^{(m)}(\mathcal{L}_T^{j}) \right) f(x) \\
& \simeq \sum_{j=1}^{n_m} \lambda_j f \left( \text{Exp}\left( \Gamma (\mathcal{L}_T^{j})\right)(x) \right)
  \end{split}
\end{equation*}
using Eq.~(\ref{eq:liepolysignaturemomentmatching}).
The error introduced when interchanging exp and $\Gamma$ 
turns out to be of the similar order with the error 
in the cubature approximation of the path level
as shown below.

Formally we define the cubature approximation at the flow level as follows.
Let $t \mapsto X_t^{x,\mathcal{L}^{j}_{\Delta}}$ for $t \in [0,1]$ 
be the deterministic process satisfying
\begin{equation}
  \label{eq:auto}
dX_{t}^{x,\mathcal{L}^{j}_\Delta} = \Gamma(\mathcal{L}^{j}_\Delta)(X_{t}^{x,\mathcal{L}^{j}_\Delta}) \,dt
\end{equation}
and 
$X_0^{x,\mathcal{L}^{j}_{\Delta}}=x$.
We define the operator 
\begin{equation}
  \label{eq:KLVflow}
\widetilde{\text{KLV}}^{(m)}
\left(\sum_{i=1}^n \kappa_i \delta_{x^i}, \Delta \right)
\equiv \sum_{i=1}^n \sum_{j=1}^{n_m} \kappa_i \lambda_j \delta_{X_{1}^{x^i, \mathcal{L}^{j}_{\Delta}}}
\end{equation}
and a sequence of discrete measure
\begin{equation*}
  \begin{split}
\widetilde{\Phi}_{\mathcal{D}}^{m,0}({\mu^0}) &={\mu^0},  \\
\widetilde{\Phi}_{\mathcal{D}}^{m,j}({\mu^0})
&=\widetilde{\text{KLV}}^{(m)}(\widetilde{\Phi}_{\mathcal{D}}^{m,j-1}({\mu^0}),s_j) 
  \end{split}
\end{equation*}
for $1 \leq j \leq k$.

Let
$\widetilde{Q}^m_{T}f(x) \equiv (\widetilde{\text{KLV}}^{(m)}(\delta_x,T),f)$
be a flow level cubature approximation,
then
the Taylor expansions of 
Eq.~(\ref{eq:nonauto})
and
Eq.~(\ref{eq:auto})
lead to
\begin{equation}
  \label{eq:flowlevelapp}
\parallel (Q^{m}_{T}-\widetilde{Q}^{m}_{T})f \parallel_\infty
\leq C \sum_{m+1 \leq \parallel I \parallel \leq 2m} T^{\parallel I \parallel/{2}} 
\parallel V_I f \parallel_\infty
\end{equation}
for a smooth $f$, where $C$ is a constant depending on $m$, $d$, 
$\mathbb{Q}^m_1$
and
$\widetilde{\mathbb{Q}}^m_1$
\cite{kusuoka2001approximation}.

\begin{theorem}
The error estimate
\begin{equation}
  \label{eq:KLVconvergeceflow}
 \sup_{x \in \mathbb{R}^N} 
 \left\lvert 
 P_Tf(x)
-({\widetilde{\Phi}^{m,k}_{\mathcal{D}(\gamma)}(\delta_x)},f) \right\rvert 
C \parallel \nabla f \parallel_\infty 
T^{1/2}
k^{-(m-1)/2}
\end{equation}
is satisfied
for a Lipschitz continuous $f$
when $\gamma>m-1$.
\end{theorem}

Eq.~(\ref{eq:KLVconvergeceflow}) 
is obtained using Eq.~(\ref{eq:flowlevelapp})
and
demonstrates that
for a suitable partition
the bounds for the approximation at flow and path level have the same rate of convergence
in view of Eq.~(\ref{eq:KLVconvergece}).
Therefore the path level operator ${\text{KLV}}^{(m)}$ 
can be replaced
by the flow level operator 
$\widetilde{\text{KLV}}^{(m)}$ 
without harming the order of accuracy.
By doing this to the PPF at the path level,
we define the 
\textit{PPF at the flow level}
by the successive algorithm that produces
$\widetilde{\pi}^{\text{PPF}}_{n|n-1}$
and
$\widetilde{\pi}^{\text{PPF}}_{n|n}$.
Furthermore,
by replacing $\text{KLV}^{(m)}$ by $\widetilde{\text{KLV}}^{(m)}$,
we define the 
\textit{APPF at the flow level}
that produces
$\widetilde{\pi}^{\text{APPF}}_{n|n-1}$ and
$\widetilde{\pi}^{\text{APPF}}_{n|n}$
instead of
${\pi}^{\text{APPF}}_{n|n-1}$ and
${\pi}^{\text{APPF}}_{n|n}$.

\subsection{Adaptive partition}
\label{sec:adapart}
Recall that
the regularity estimate of
Eq.~(\ref{eq:ugf})
for the $\text{Lip}(1)$ function 
is used to obtain
the higher order error bound
of Eq.~(\ref{eq:KLVconvergece}).
It implies that the Kusuoka partition in Eq.~(\ref{eq:Kpartition})
is suitable to accurately integrate 
any $\text{Lip}(1)$ function.

For a given test function $f$ that
might have different regularity
from the Lipschitz continuity,
one can make use of the heat kernel
$P_t$ 
to evolve cloud of particles
so that one step error is
within a given degree of accuracy, i.e.,
\begin{equation}
  \begin{split}
	\label{eq:klverror}
  \parallel (P_{s_j}-Q^m_{s_j})P_{T-t_j}f \parallel_\infty  < \epsilon
  \end{split}
\end{equation}
for some $\epsilon>0$.
We define 
an adaptive partition
$\mathcal{D}({\epsilon},f) = \{t_j\}_{j=0}^k$
to be a time discretization for which
each
$s_j=t_j-t_{j-1}$ is the supremum among
the ones
satisfying
Eq.~(\ref{eq:klverror}).
Because
$P_tf$ becomes smoother as $t$ increases,
the sequence $\{s_j\}_{j=1}^k$
tends
to decrease monotonically,
i.e.,
$s_1 \geq s_2 \geq \cdots \geq s_k$.
The upper bound of the total error due to the adaptive partition
is given by
\begin{equation*}
\lvert P_Tf(x) -(\Phi^{m,k}_{\mathcal{D}}(\delta_x),f) 
\rvert
< k \epsilon
\end{equation*}
from Eq.~(\ref{eq:inequality1}).

In order to find
$\mathcal{D}({\epsilon},f)$,
the smooth function $P_{T-t_j}f$ has to be approximated.
One solution is to compute 
$P_{T-t_j}f(z_i)$ for finitely many $z_i \in \mathbb{R}^N$
(which can be done by using the KLV method or the Monte-Carlo method)
and 
to use the interpolation scheme developed in the scattered data approximation
\cite{wendland2005scattered}.

When $f$ is a Lipschitz function,
adaptive partition can be analyzed by the Lipschitz norm.
Note that the inequalities
	\begin{align}
  \label{eq:adapartlip}
\parallel (P_{s_j}-Q^m_{s_j})P_{T-t_j}f \parallel_\infty  
& 
 \leq 
C s_j^{\frac{m+1}{2}} \sup_{\parallel I \parallel = m+1,m+2}  \parallel
V_IP_{T-t_j}f\parallel_\infty \nonumber \\
& 
\leq C s_j^{\frac{m+1}{2}} \parallel P_{T-t_j}f\parallel_{\text{Lip}(m+2)}
	\end{align}
hold
from Eq.~(\ref{eq:errorbound})
and from
the definition of the Lipschitz norm.
Suppose that we have the axiom
\begin{equation}
  \label{eq:Lip}
\parallel P_tf \parallel_{\text{Lip}(\rho')}
\leq \frac{K}{t^{\alpha}}
\parallel f \parallel_{\text{Lip}(\rho)}
\end{equation}
where $K$ and $\alpha$
are determined by $\rho$ and $\rho'$,
which is a generalisation
of the regularity estimate of 
Eq.~(\ref{eq:ugf})
or
Eq.~(\ref{eq:uh}).
One can use 
Eq.~(\ref{eq:adapartlip})
and
Eq.~(\ref{eq:Lip}) 
to quantify $s_j$ 
satisfying Eq.~(\ref{eq:klverror})
in terms of 
the Lipschitz norm of $f$.

\subsection{Adaptive recombination}
\label{sec:adarecomb}
Similarly with Eq.~(\ref{eq:klverror}),
we consider the condition
\begin{equation}
	\label{eq:adarecomb}
  \left\lvert  
  \left({\Phi^{m,j}_{\mathcal{D},(u,r)}} (\mu^0)
  -\widehat{\Phi}^{m,j}_{\mathcal{D},(u,r)} (\mu^0)
,P_{T-t_j}f \right)
\right\rvert < \theta
\end{equation}
where $\mu^0= \delta_x$ for the PDE problem
or $\mu^0 = \pi^{\text{PPF}}_{n-1|n-1}$ for the PPF,
given some $\theta>0$.
We define the adaptive recombination 
by the algorithm that uses the maximum value of $u$,
for fixed recombination degree $r$,
among the ones
satisfying Eq.~(\ref{eq:adarecomb}).
Combining 
Eq.~(\ref{eq:klverror})
and
Eq.~(\ref{eq:adarecomb})
yields the error bound
\begin{equation*}
 \lvert P_Tf(x) -(\Phi^{m,k}_{\mathcal{D},(u,r)}(\delta_x),f) 
 \rvert 
 < k(\epsilon+\theta)
\end{equation*}
from Eq.~(\ref{eq:errorineq}).

The implementation 
of the adaptive recombination
does not require to specify the size of patches:
it is enough to keep making
the size of patches smaller until
Eq.~(\ref{eq:adarecomb}) is satisfied.
The adaptive recombination is useful,
particularly in high dimensions,
together with the Morton ordering
\cite{morton1966computer}.
Instead of balls, the method uses boxes
as container of the patched particles.

Assume that the particles are in $[0.5, 1)^N$, a box of $N$ dimension.
In double-precision floating-point format,
any $z^i \in [0.5, 1)$ is expressed in terms of $\{b^i_j\}_{j=1}^{52}$
where $b^i_j$ is either $0$ or $1$
in a way that $z^i = (1/2)\times(1+\sum_{j=1}^{52}b^i_j2^{-j})$.
Then the point $(z^1,\cdots,z^N)$ in $N$-dimension can be expressed by $52 \times N$ binary numbers.
Interleaving the binary coordinate values yields binary values. 
Connecting the binary values in their numerical order produces 
the Morton ordering.
Then an appropriate coarse-graining leads to the subdivision of a box.
For examples, when $N=2$,
the binary value corresponding $(z^1,z^2)$ is $b^1_1b^2_1b^1_2b^2_2 \cdots b^1_{52}b^2_{52}$.
The point is 
in first quadrant if $(b^1_1,b^2_1)=(1,1)$,
in second quadrant if $(b^1_1,b^2_1)=(0,1)$,
in third quadrant if $(b^1_1,b^2_1)=(0,0)$
and
in fourth quadrant if $(b^1_1,b^2_1)=(1,0)$.
Applying this classification to a number of particles 
produces $2^2$ disjoint subsets of clustered particles.
Similarly, using $b^1_1b^2_1b^1_2b^2_2$ and ignoring the rest subgrid scales
gives $4^2$ subsets when $N=2$.
Using affine transformation when the particles are not in $[0.5, 1)^N$,
the Morton ordering within floating-point context provides an efficient way to patch the particles.

Note that the inequality
\begin{equation}
  \label{eq:adarecomblip}
\left\lvert  
  \left({\Phi^{m,j}_{\mathcal{D},(u,r)}} (\mu^0)
  -\widehat{\Phi}^{m,j}_{\mathcal{D},(u,r)} (\mu^0) ,P_{T-t_j}f \right)  \right\rvert 
\leq \frac{Cu_j^{r_j+1}}{(r_j+1)!}\parallel P_{T-t_j}f\parallel_{\text{Lip}(r_j+1) } 
\end{equation}
holds
from
Eq.~(\ref{eq:Liperror})
when $f$ is a Lipschitz function.
One can use 
Eq.~(\ref{eq:adarecomblip})
and
Eq.~(\ref{eq:Lip})
to
quantify
the patch size $u_j$ 
satisfying 
Eq.~(\ref{eq:adarecomb})
for a fixed $r_j$
in terms of
the Lipschitz norm of $f$.

\section{Numerical simulations}
\label{sec:numerical}
We here perform numerical simulations 
to examine the efficiency and accuracy of the proposed filtering approaches.
We introduce the test model in subsection~\ref{sec:tmodel}
and 
apply the Kalman filter to obtain its exact solution in subsection~\ref{sec:kfilter}.
Subsequently 
we obtain particle approximations
from a Monte-Carlo sampling in subsection~\ref{sec:montecarlosamples}
and 
from
PPF and APPF with cubature on Wiener space of degree $m=5$ in subsection~\ref{sec:ppfappfd5}.
Finally, 
in subsection~\ref{sec:ppfappfd7},
we study the prospective performance of 
PPF and APPF with cubature on Wiener space of degree $m=7$.

\subsection{Test model}
\label{sec:tmodel}
Consider the 
stochastic differential equation
	\begin{equation}
	  \label{eq:L63}
dX=
d\left( 
\begin{array}{c}
x^1   \\
x^2  \\
x^3
\end{array} 
\right)
=  
  \left[-\Lambda X+ 
  a_0
\left( 
\begin{array}{c}
0   \\
-x^1x^3  \\
x^1x^2
\end{array} 
\right)
\right]dt+g I_3
dW
	\end{equation}
	in three dimension
where
$\Lambda=
	\left( 
	\begin{array}{ccc}
	\sigma   & -\sigma  & 0 \\
	-\rho  & 1 & 0\\
	0 & 0 & \beta
	\end{array} 
	\right)
$,
$dW= (dW_1\; dW_2\; dW_3)^t $
and $I_3$ denotes the $3 \times 3$ identity matrix.
Here the superscript $t$ denotes the transpose.
When $a_0=1$, Eq.~(\ref{eq:L63}) is the Lorenz-$63$ model
that has been intensively studied in the data assimilation community
\cite{lorenz1963deterministic,
miller1999data,
xiong2006note}.
When $a_0=0$, 
Eq.~(\ref{eq:L63}) is
the Ornstein-Uhlenbeck process
\cite{uhlenbeck1930theory}
and this linear process is our test model.
The parameter values $\sigma=1$, $\rho=0.28$, $\beta=8/3$
and $g=0.5$
are used.

For the observation process,
we take
the identity function
for $\varphi(\cdot)$ in Eq.~(\ref{eq:observation}),
and 
\begin{equation}
\label{eq:linearobservation}
Y_n  = X_n + \eta_n, \quad \eta_n \sim \mathcal{N}(0,R_n).
\end{equation}
The inter-observation time is $T=0.5$ and 
the cases 
are studied
in which
the covariance of observation noise is $R_n=R \times I_{3}$
for a number of values $R = 10^{-1}, 10^{-2}$ and $10^{-3}$.

\subsection{Kalman filter}
\label{sec:kfilter}
For Eq.~(\ref{eq:L63}) where $a_0=0$ and Eq.~(\ref{eq:linearobservation}),
the conditioned measure is Gaussian 
and
$\pi_{n|n'}=\mathcal{N}(M_{n|n'},C_{n|n'}) $ 
can be obtained using the Kalman filter.
In this case, the prior covariance $C_{n|n-1}$ satisfies 
the Riccati difference equation and its solution converges
as $n$ increases
(see \cite{bitmead1991riccati} for the conditions).
We take the measure $\pi_{0|0}$ so that
$C_{n|n-1}$ and $C_{n|n}$ do not depend on $n$.
We see that
the diagonal element of $C_{n|n-1}$ are about $10^{-1}$ for all cases of $R=10^{-1},10^{-2},10^{-3}$.
The diagonal element of $C_{n|n}$ are about $10^{-1}$ when $R=10^{-1}$,
about $10^{-2}$ when $R=10^{-2}$ and about $10^{-3}$ when $R=10^{-3}$.

\begin{figure}
\centerline{
\includegraphics[width=0.34\textwidth]{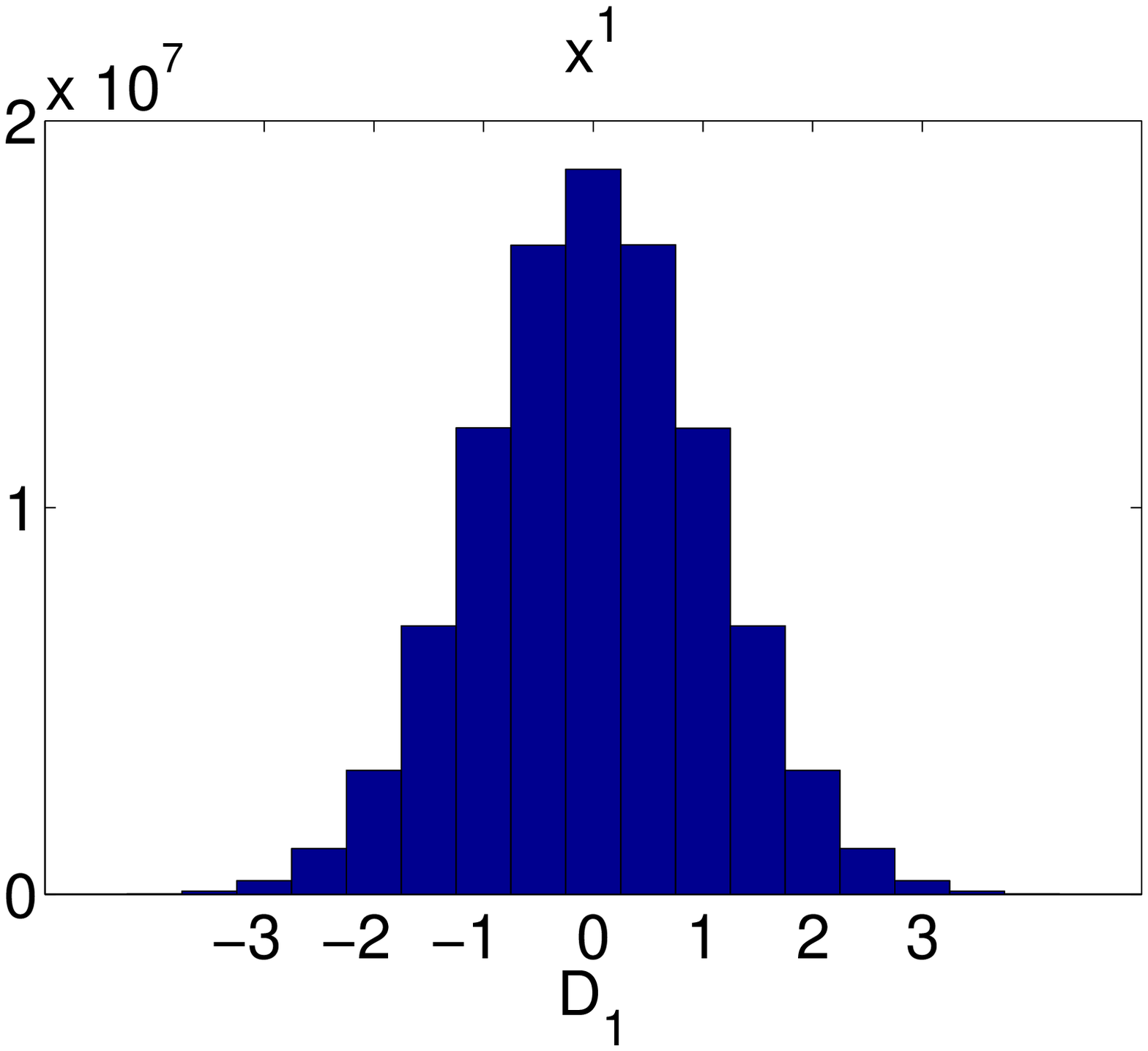}
\includegraphics[width=0.34\textwidth]{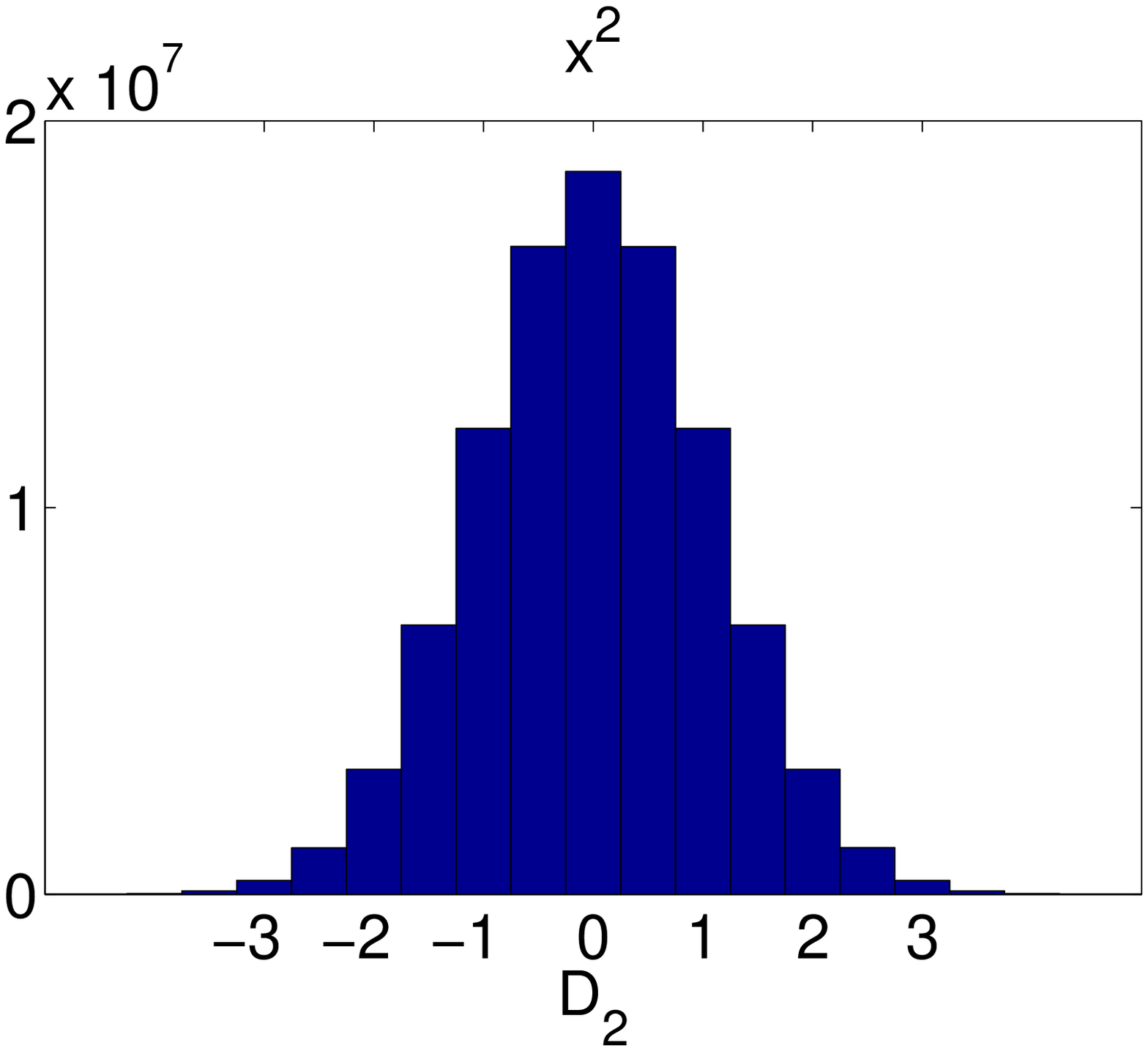}
\includegraphics[width=0.34\textwidth]{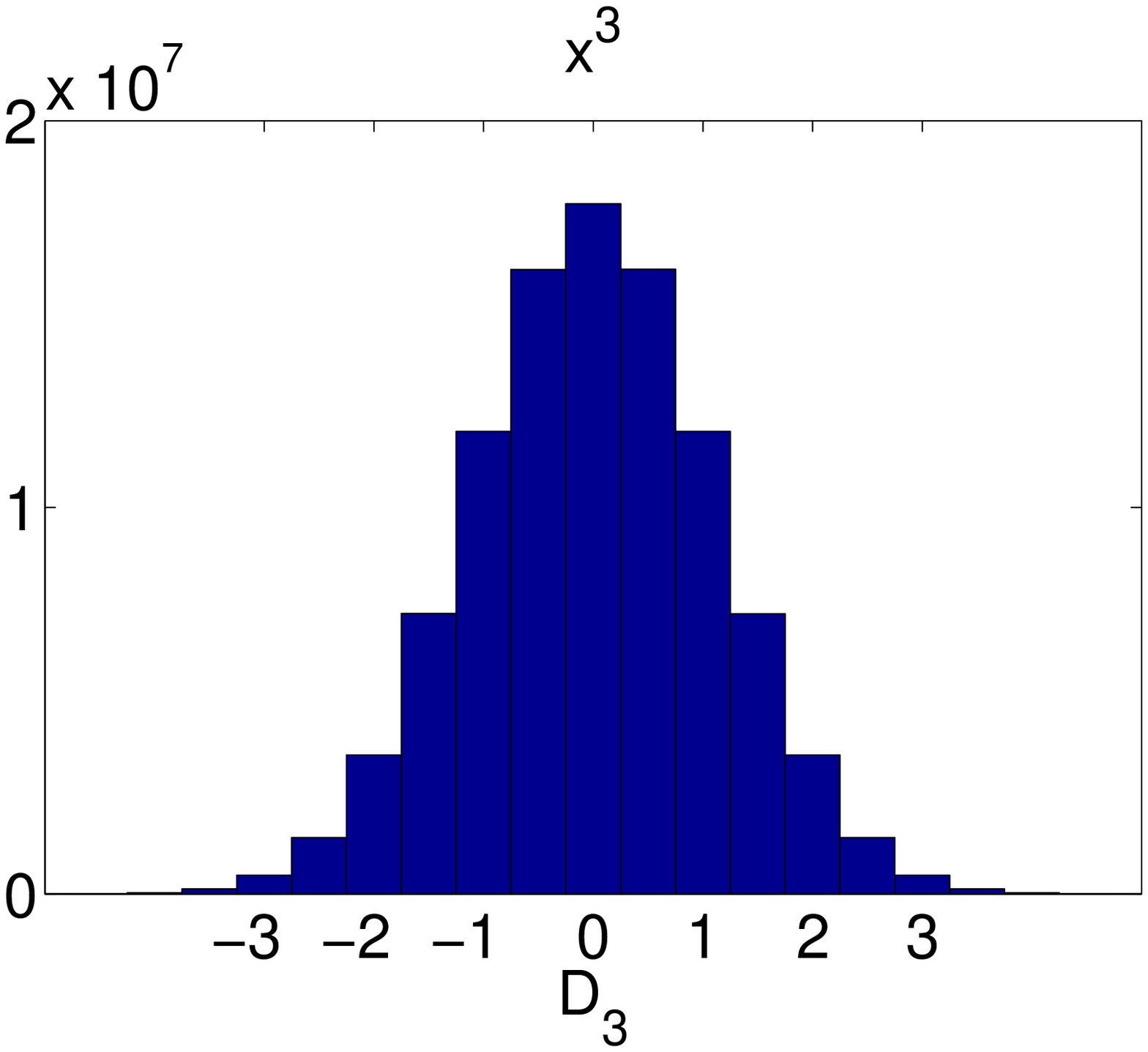}
}
\caption{
 The distribution of normalized distances between the observation and the prior mean
 when the noise covariance is $R_n = 10^{-2} \times I_3$.
}
\label{fig:observations} 
\end{figure}

We apply the Kalman filter for $ 1\leq n \leq 10^8$ 
and calculate 
the values of $D_1, D_2$ and $D_3$ satisfying
$y_n={M}_{n|n-1}+(D_1\; D_2\; D_3)^t \cdot \sqrt{\text{diag}({C}_{n|n-1})}$.
The histograms in 
Fig.~\ref{fig:observations}
show the distribution of
these normalized distances between the observation and the prior mean
when $R=10^{-2}$
(the cases of $R=10^{-1}$ and $R= 10^{-3}$ are similar and not shown).
One can see that
most of the observations are within two times the standard deviations from the prior
mean in each coordinate.
Among the cases of $10^8$,
there are $4,592,208$ cases for which $|D_i|>1$ for all $i=1,2,3$ at the same time.
There are $37,574$ cases for which $|D_i|>2$ for all $i$ at the same time,
and $60$ cases for which $|D_i|>3$ for all $i$ at the same time.
In the following, we study the cases when the parameter value of $D \equiv D_1=D_2=D_3$ 
is $1$, $2$ and $3$. These three cases correspond to normal, exceptional and rare event,
respectively.

\subsection{Monte-Carlo samples}
\label{sec:montecarlosamples}
In order to quantify the accuracy of the discrete measures
given in the form of Eq.~(\ref{eq:weighteddiscretemeasure}),
we define the $L^2$ norm of the moment as the following.
Let $\mathcal{C}^p$ be the $p$-th moment of $X=(x^1 \; x^2 \; x^3)^t$, i.e.,
\begin{equation*}
  \mathcal{C}^p_{i_1, \cdots, i_p} = 
  \mathbb{E}
  \left(
  \prod_{j=1}^p \left(x^{i_j}
  -\mathbb{E}( {x}^{i_j} )
  \right) 
  \right)
\end{equation*}
where $i_j = 1,2, 3$.
The $L^2$ norm of $\mathcal{C}^p$
is defined by
\begin{equation}
  \label{eq:momentnorm}
  \parallel \mathcal{C}^p \parallel_2
  \equiv
	\left(
	\sum_{i_1, \cdots, i_p =1 }^3 | \mathcal{C}^p_{i_1, \cdots, i_p}|^2 
	\right)^{1/2}.
\end{equation}
When $p=1$, 
Eq.~(\ref{eq:momentnorm}) is the Euclidean norm of the vector.
It is the Frobenius norm of the matrix when $p=2$.
The relative root mean square error 
\begin{equation}
  \label{eq:rmse}
\text{rmse} \%
\equiv
\parallel \mathcal{C}^p - \widehat{\mathcal{C}}^p \parallel_2/ \parallel \mathcal{C}^p \parallel_2
\end{equation}
is used to measure the accuracy of the moment approximations,
where $\widehat{\mathcal{C}}^p$ is the $p$-th moment of a particle approximation.

We perform a Gaussian sampling from the prior and the posterior
(the samples are not obtained from the integration of Eq.~(\ref{eq:L63})
but drawn from the analytic solution).
We also apply bootstrap reweigthing to the prior samples 
so that it approximates the posterior.
The values of 
Eq.~(\ref{eq:rmse})
are depicted in
Figs.~\ref{fig:Ra},
\ref{fig:Rc},
\ref{fig:Re},
\ref{fig:Rg}
when $R=10^{-2}$, $D=1,2,3$
and
in
Figs.~\ref{fig:Da},
\ref{fig:Db},
\ref{fig:Dd},
\ref{fig:De},
\ref{fig:Dg},
\ref{fig:Dh}
when $D=1$,
$R=10^{-1},10^{-2},10^{-3}$.
As the number of samples $M$ increases, the error of
Eq.~(\ref{eq:rmse})
asymptotically
behaves $M^{-1/2}$
in all cases.
The results will be used as the benchmark calculations
for the accuracy of PPF and APPF.

\subsection{PPF and APPF with cubature on Wiener space of degree $5$}
\label{sec:ppfappfd5}
We here apply PPF and APPF at the flow level.
In case of $d=3$,
i.e.,
the system is driven by three independent white noises,
cubature on Wiener space of degree $m=3$ and $m=5$,
with support size $n_m=6$ and $n_m=28$ respectively, are available.
We use the KLV operator with degree $m=5$.

We use the adaptive partition for both PPF and APPF.
In order to do that, the partition
$\mathcal{D}(\epsilon,g^{y_n})$
satisfying Eq.~(\ref{eq:klverror})
with $\widetilde{Q}^m_{s_j}$ in place of ${Q}^m_{s_j}$
is analytically obtained
for the system of Eq.~(\ref{eq:L63}) where $a_0=0$.
Note that the likelihood 
$g^{y_n}$ is the density function of $\mathcal{N}(y_n,R_n)$
and that the adaptive partition does not depend on $y_n$ but on $R_n$
or Lipschitz norm of $g^{y_n}$.
The number of iterations $k$ as a function of $\epsilon$ and $R$ is listed in
table~\ref{tab:partm5}.

\begin{table}
  \renewcommand{\arraystretch}{1.3} 
  \caption{The number of adaptive partition $k$ for KLV with $m=5$} 
  \label{tab:partm5} 
  \centering \begin{tabular}{ c | c  c cc} 
	& $\epsilon=10^{-2}$ & $\epsilon=10^{-3}$ & $\epsilon=10^{-4}$ & $\epsilon=10^{-5}$ \\ \hline
	$R=10^{-1}$ & 7 & 31 & 102 & 344  \\ 
	$R=10^{-2}$ & 10 & 29 & 101 & 330  \\ 
	$R=10^{-3}$ & 20 & 48 & 120 &  329 \\ 
  \end{tabular} 
\end{table}

For the recombination of the PPF,
we use a variant of the adaptive recombination.
It requires to satisfy 
Eq.~(\ref{eq:adarecomb})
with $f=g^{y_n}$ for all $y_n \in \mathbb{R}^N$, i.e.,
\begin{equation}
	\label{eq:supadarecomb}
	\sup_{y_n}  \left\lvert  
  \left({\Phi^{m,j}_{\mathcal{D},(u,r)}} (\mu^0)
  -\widehat{\Phi}^{m,j}_{\mathcal{D},(u,r)} (\mu^0)
  ,P_{T-t_j}g^{y_n} \right)
\right\rvert < \theta
\end{equation}
so that the recombination does not depend on $y_n$
but on $R_n$.
We choose the recombination degree $r=5$
and simulate the PPF for the cases of $\epsilon=10^{-2}, 10^{-3}$
with $\theta = 0.3 \times \epsilon$.

For the APPF, 
the tolerance $\tau$ has to be specified
in addition to the parameters $\{ \epsilon, \theta\}$.
The value of $\tau$ varies in each case,
but we choose it so that
the ADA operator 
in Eq.~(\ref{eq:adaptivesequenceofdiscretemeasure})
allows $1/4 \sim 1/3$ part of particles leap to the next observation time
for all iterations
except the first and last few steps.
The remaining particles are reduced by the adaptive recombination,
i.e.,
the recombination satisfies
\begin{equation}
	\label{eq:appfadarecomb}
  \left\lvert  
  \left({\Phi^{m,j}_{\mathcal{D},(u,r),\tau}} (\mu^0)
  -\widehat{\Phi}^{m,j}_{\mathcal{D},(u,r),\tau} (\mu^0)
  ,P_{T-t_j} g^{y_n} \right)
\right\rvert < \theta
\end{equation}
where $\mu^0 = \widetilde{\pi}^{\text{APPF}}_{n-1|n-1}$.
We again choose the recombination degree $r=5$
and simulate the APPF for the cases of $\epsilon=10^{-2}, 10^{-3}$
with $\theta = 0.3 \times \epsilon$.

\begin{figure}
\centerline{
  \subfigure[PPF using adaptive partition with $\epsilon=10^{-2}$]
  {\label{fig:Tppf}\includegraphics[width=0.49\textwidth]{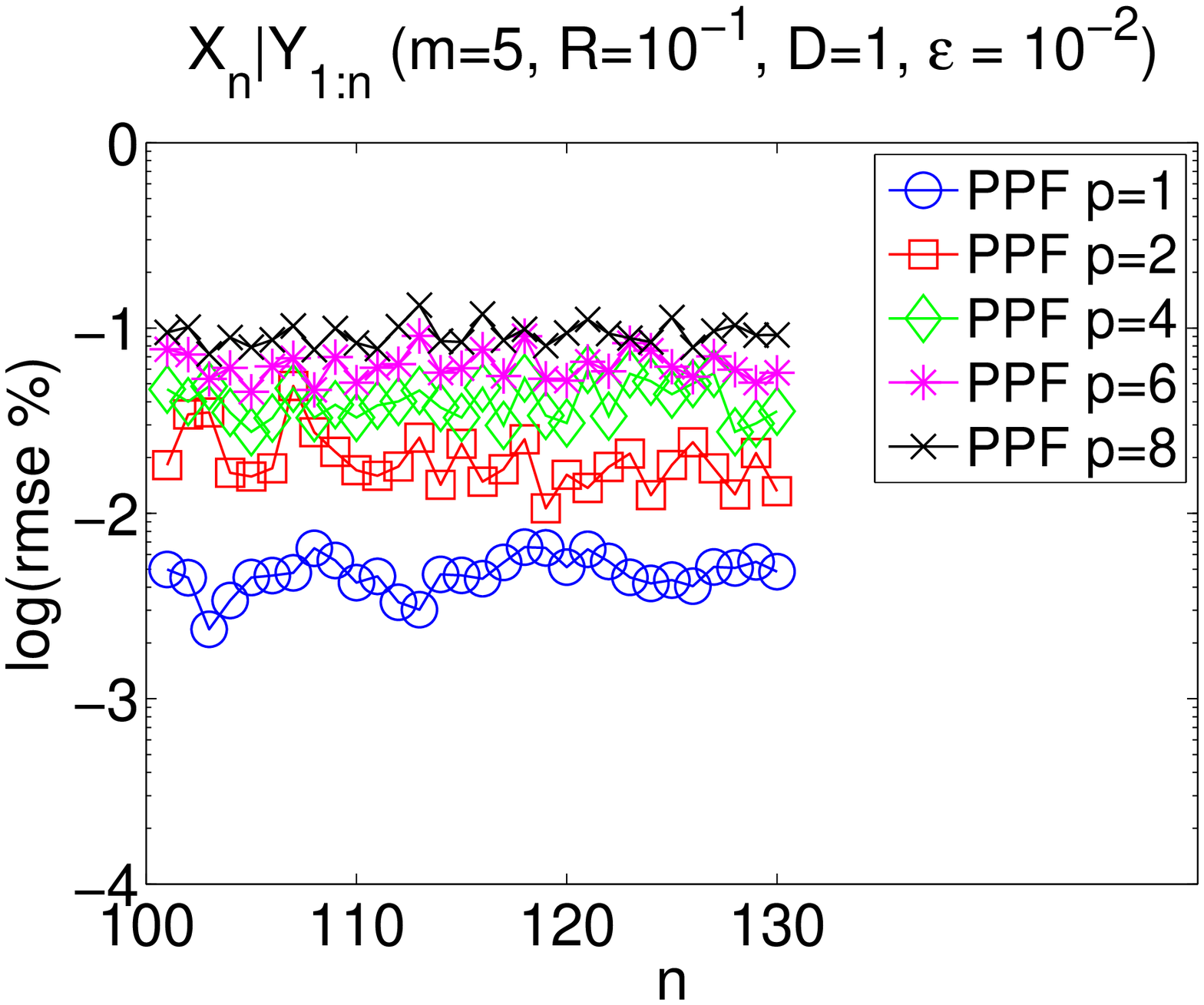}}
\subfigure[APPF using adaptive partition with $\epsilon=10^{-3}$]
{\label{fig:Tappf}\includegraphics[width=0.49\textwidth]{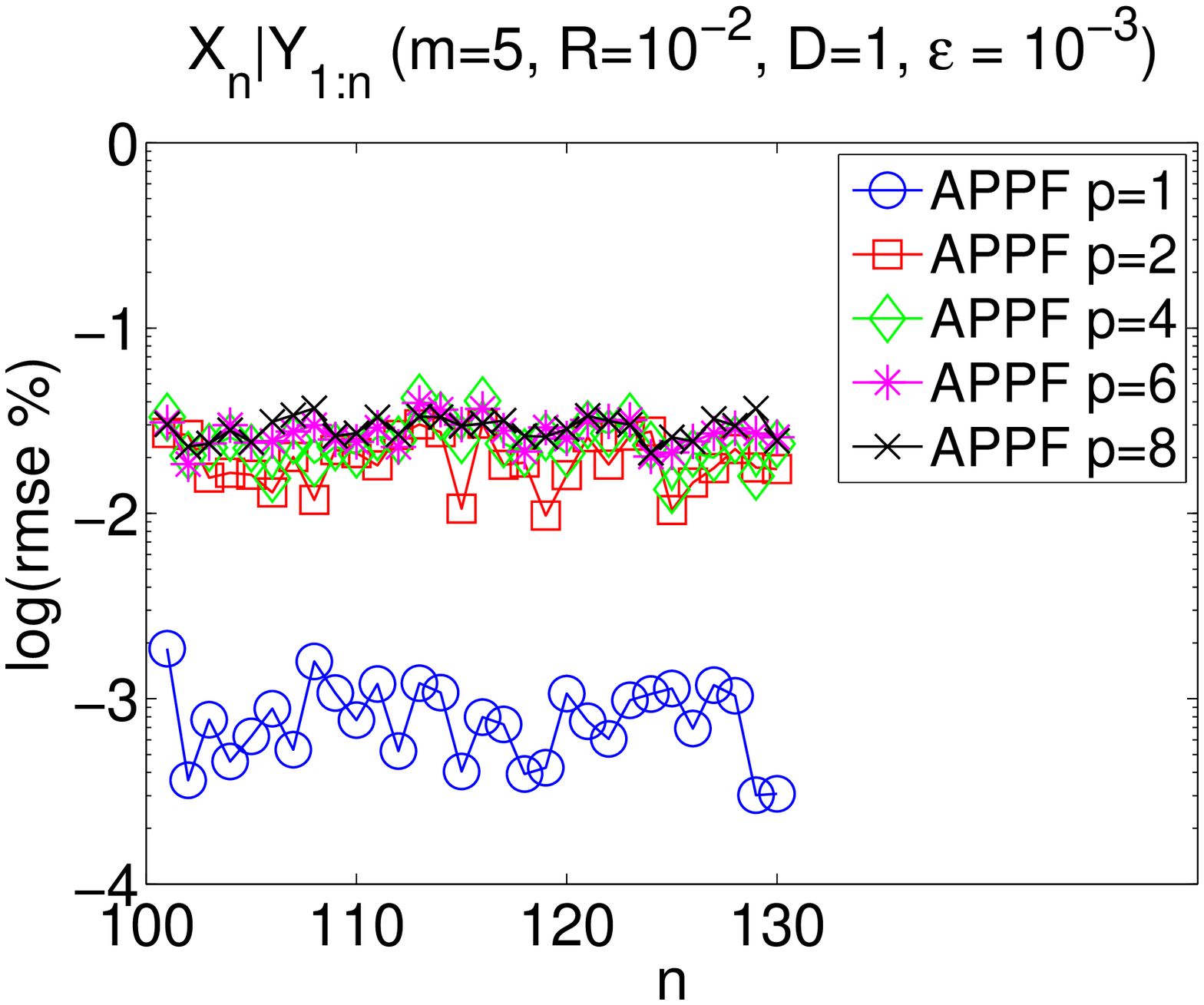}}
}
\caption{
The relative $L^2$ errors for the $p$-th moments of
the evolutionary posterior.
}
\label{fig:T}
\end{figure}

With the value of $D$ being fixed,
we apply PPF and APPF to obtain 
the values of Eq.~(\ref{eq:rmse})
for the evolutionary posterior.
Fig.~\ref{fig:T} shows that
two filtering algorithms are stable and 
that our numerical error estimates can be trusted
(the rest cases produce similar plots and are not shown).

\begin{figure*}[!htp] 
\centerline{
\subfigure[unweighted prior samples]
{\includegraphics[width=0.43\textwidth]{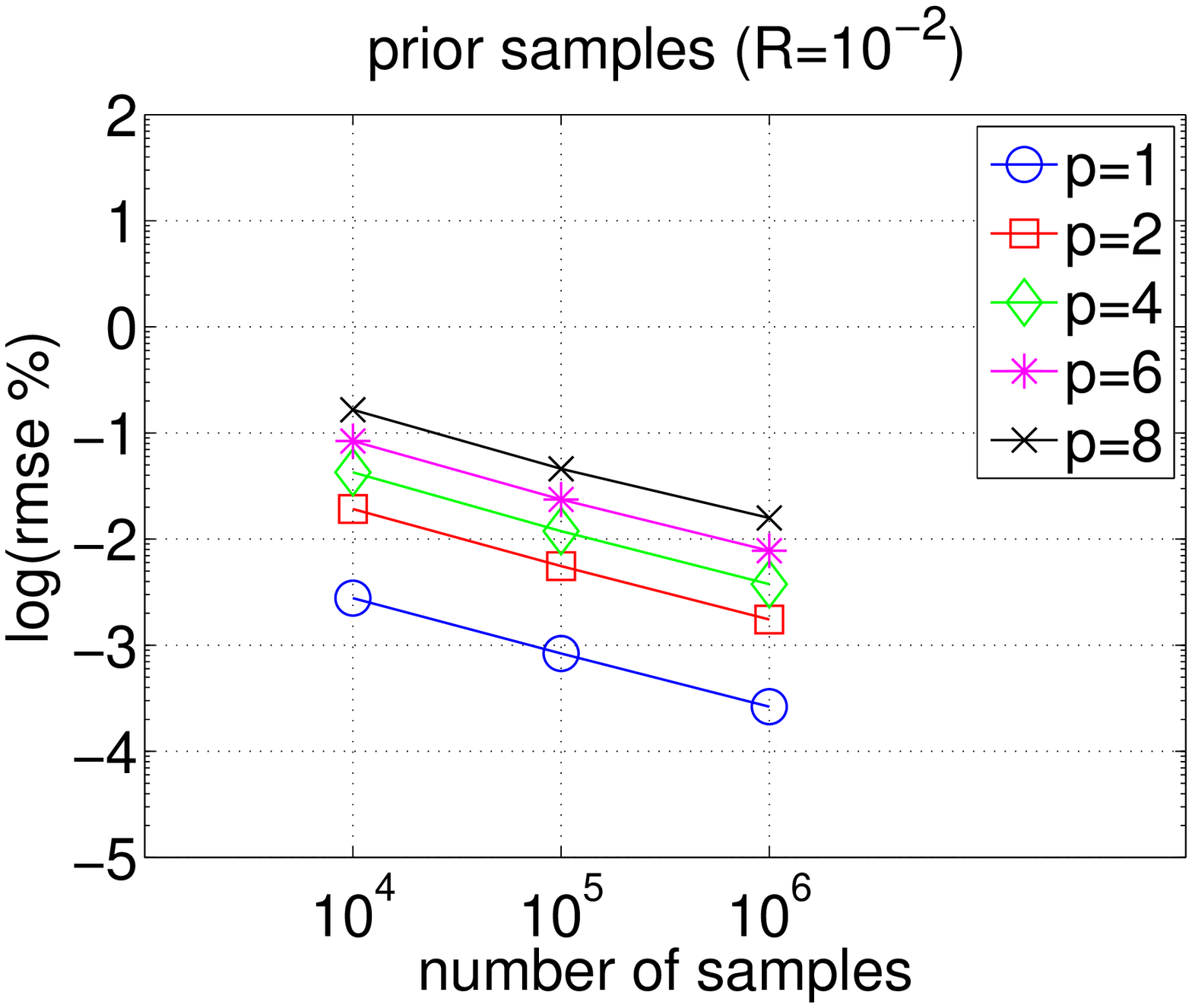} \label{fig:Ra} } 
\quad
\subfigure[cubature approximation of prior]
{\includegraphics[width=0.43\textwidth]{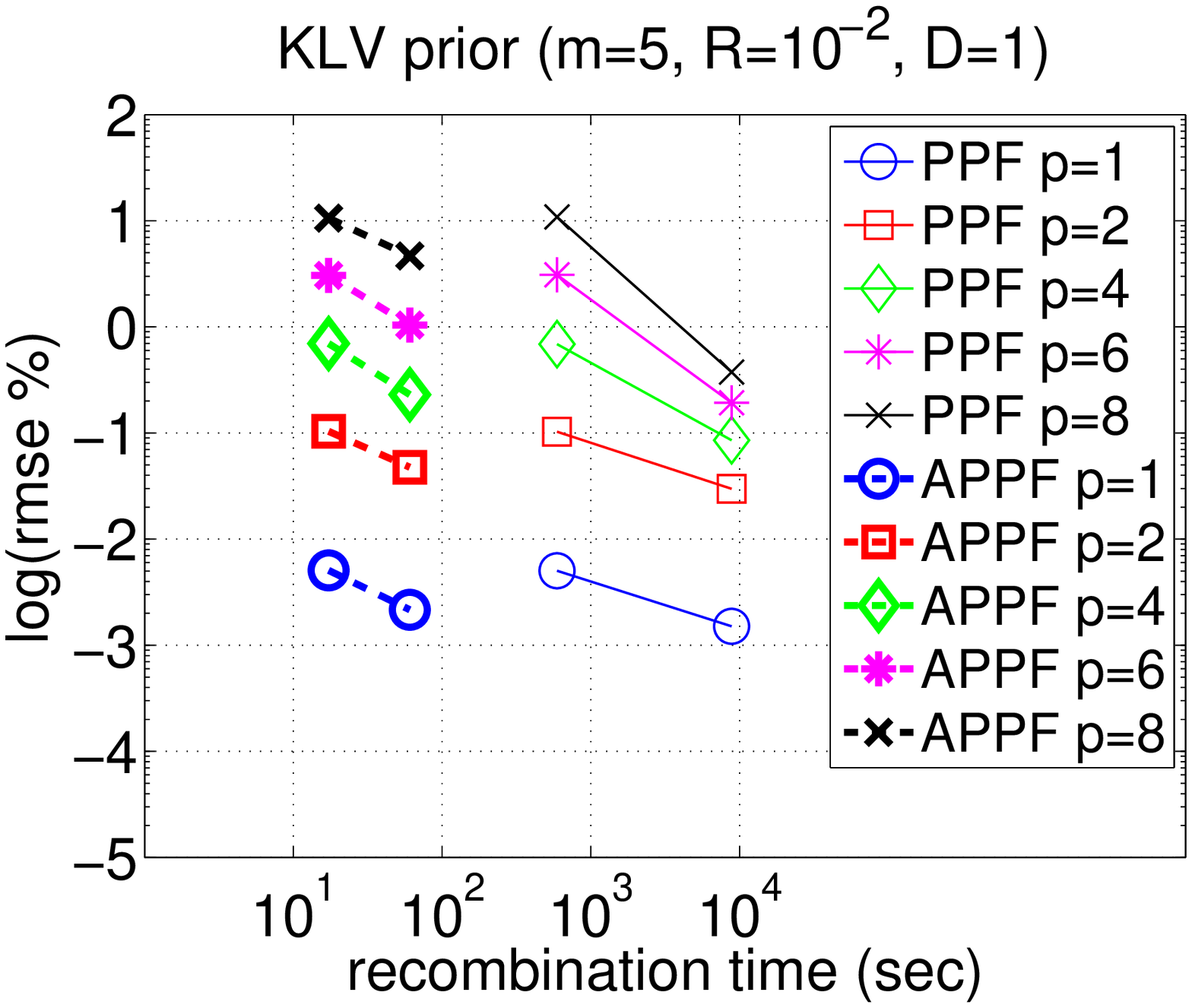} \label{fig:Rb}} 
} 
\centerline{
\subfigure[bootstrap reweighted prior samples]
{\includegraphics[width=0.43\textwidth]{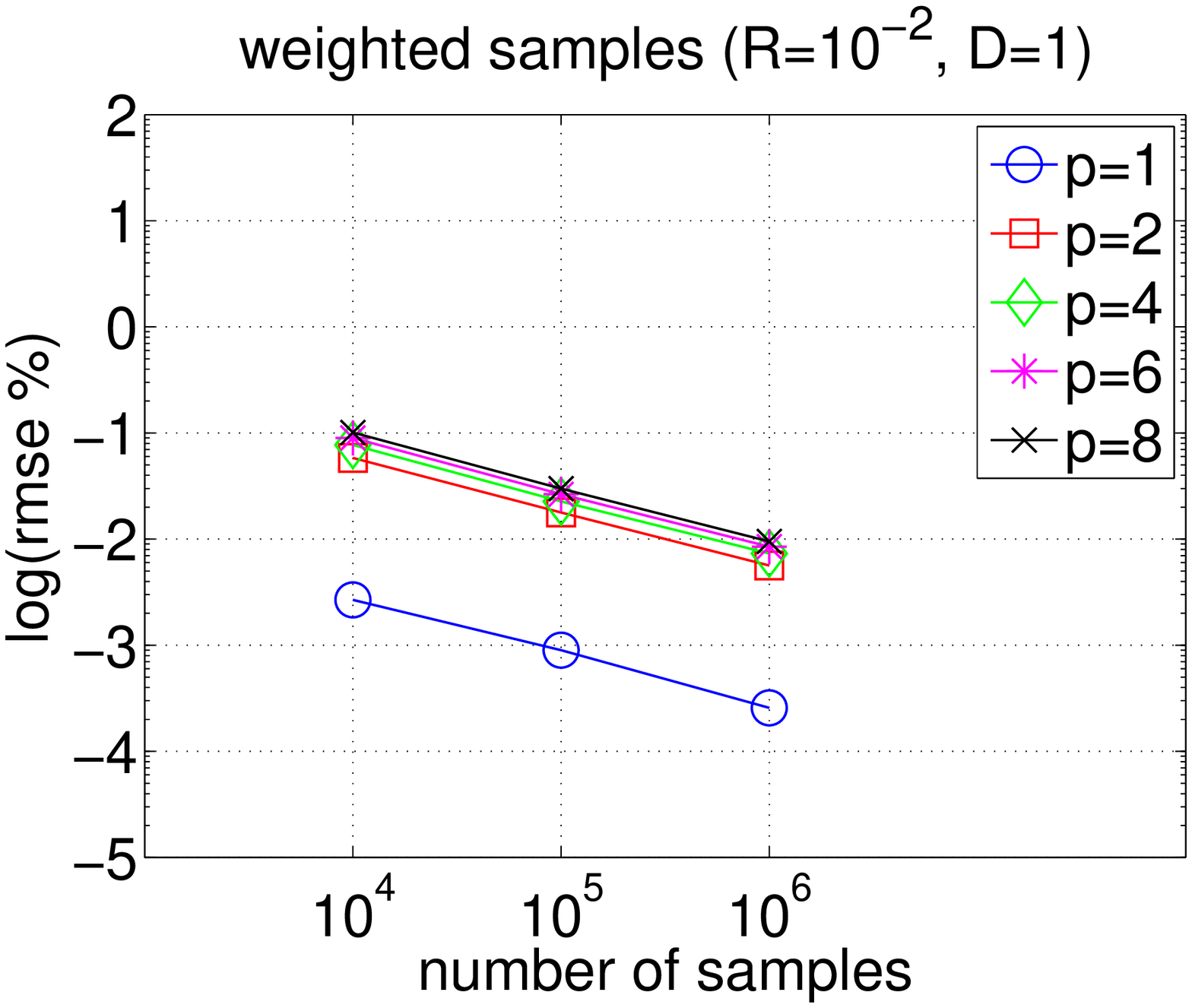} \label{fig:Rc}} 
\quad
\subfigure[cubature approximation of posterior]
{\includegraphics[width=0.43\textwidth]{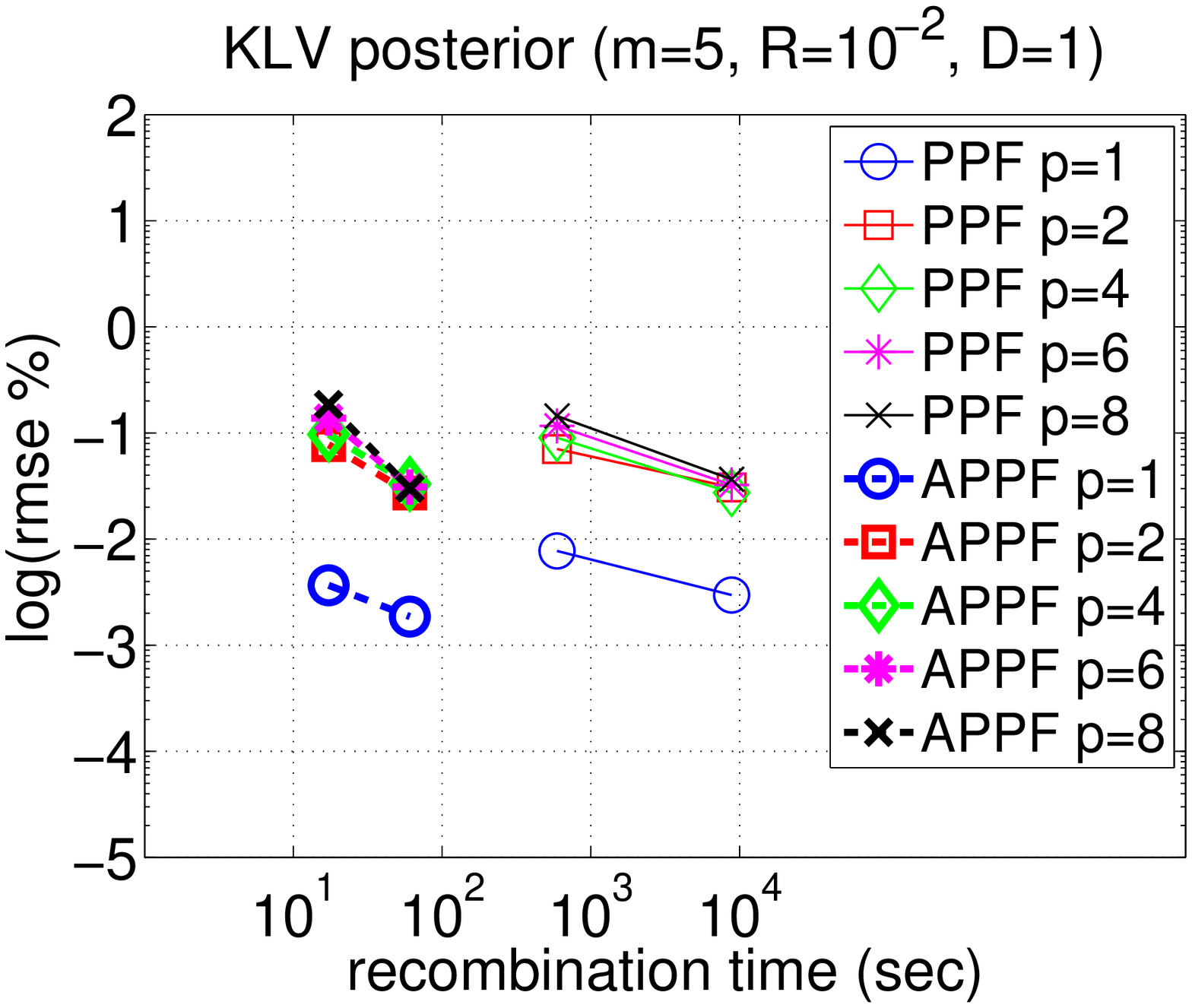} \label{fig:Rd}} 
}
\centerline{
\subfigure[bootstrap reweighted prior samples]
{\includegraphics[width=0.43\textwidth]{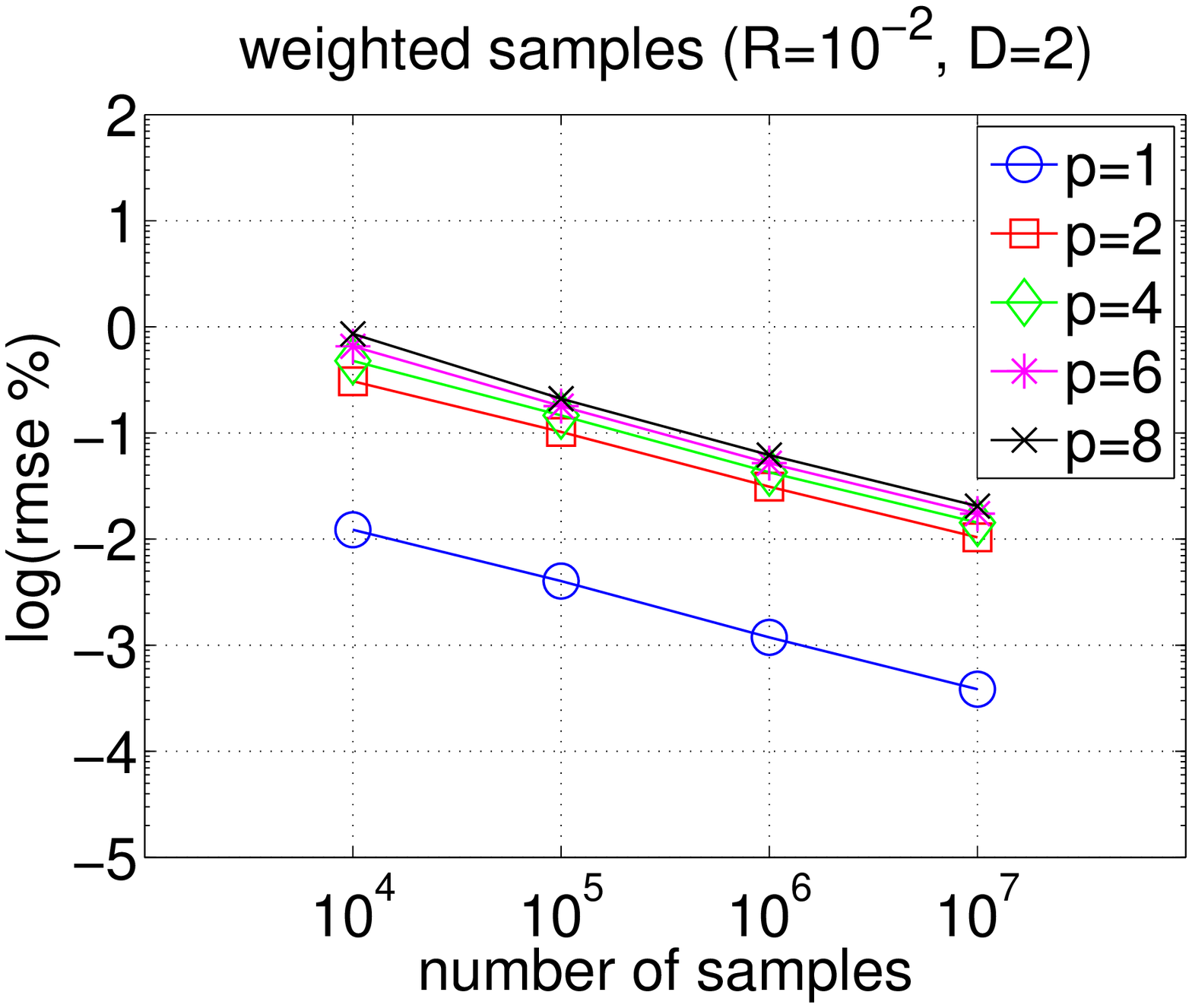} \label{fig:Re}} 
\quad
\subfigure[cubature approximation of posterior]
{\includegraphics[width=0.43\textwidth]{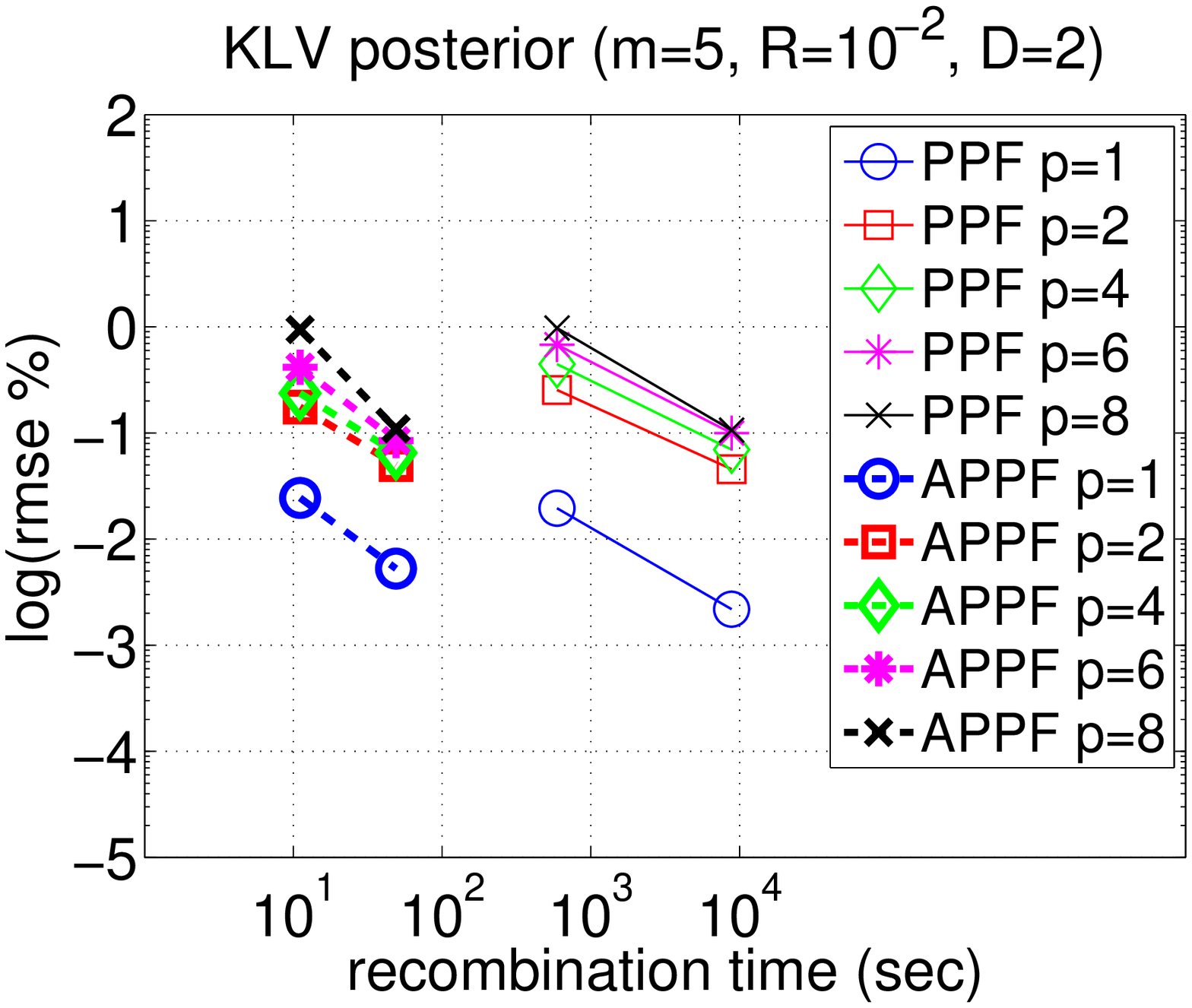} \label{fig:Rf}} 
}
\centerline{
\subfigure[bootstrap reweighted prior samples]
{\includegraphics[width=0.43\textwidth]{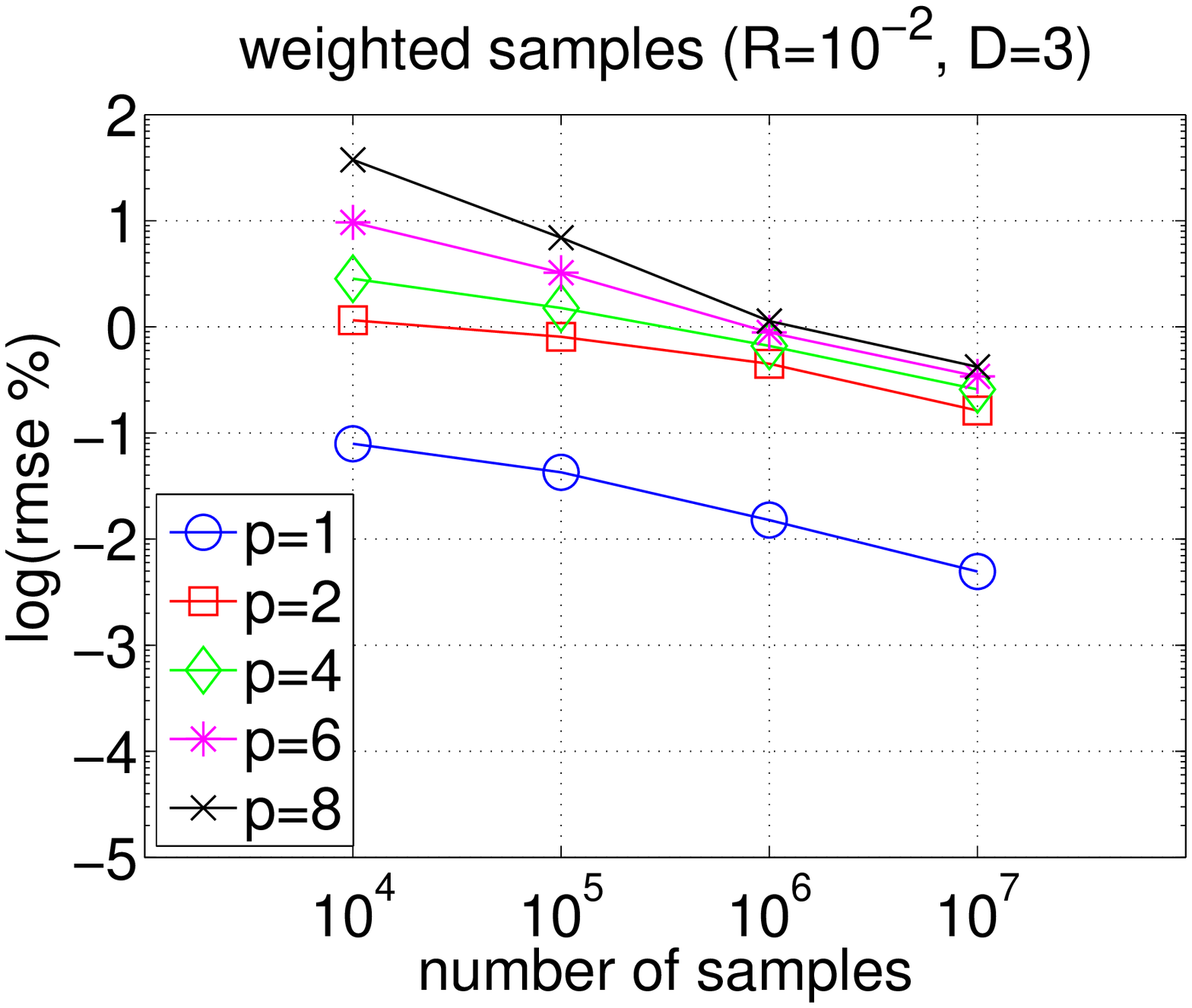} \label{fig:Rg}} 
\quad
\subfigure[cubature approximation of posterior]
{\includegraphics[width=0.43\textwidth]{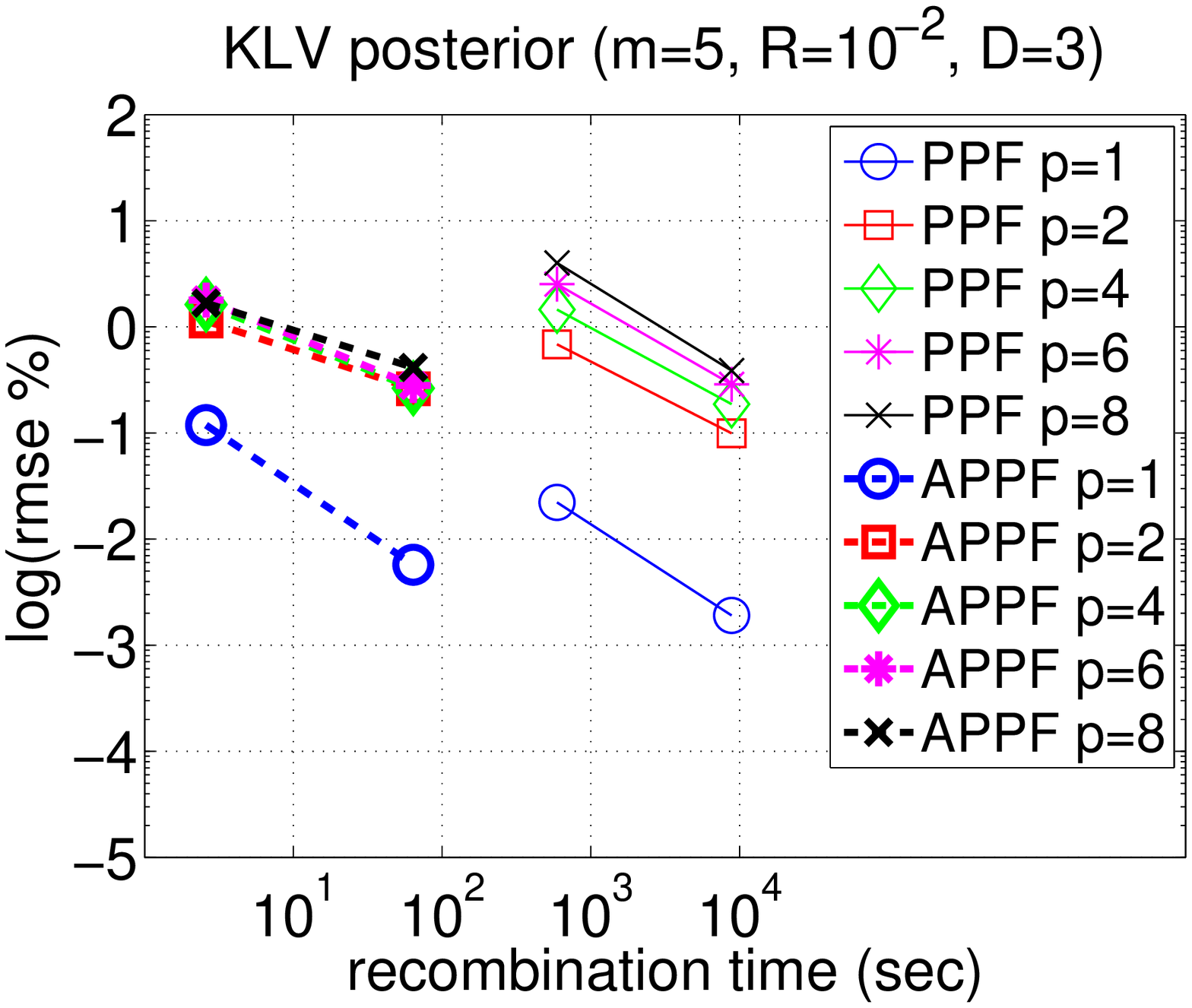} \label{fig:Rh}} 
}
\caption{
The prior and posterior approximations 
when $R=10^{-2}$ is fixed and $D=1,2,3$ varies. 
The left column is from Monte-Carlo samples and the right column is from cubature 
approximation when $\epsilon=10^{-2}, 10^{-3}$.
The top first row is for the prior and the bottom three rows are for the posterior.
} 
\label{fig:R} 
\end{figure*}

\begin{figure}
\centerline{
\subfigure[unweighted posterior samples]
{\includegraphics[width=0.39\textwidth]{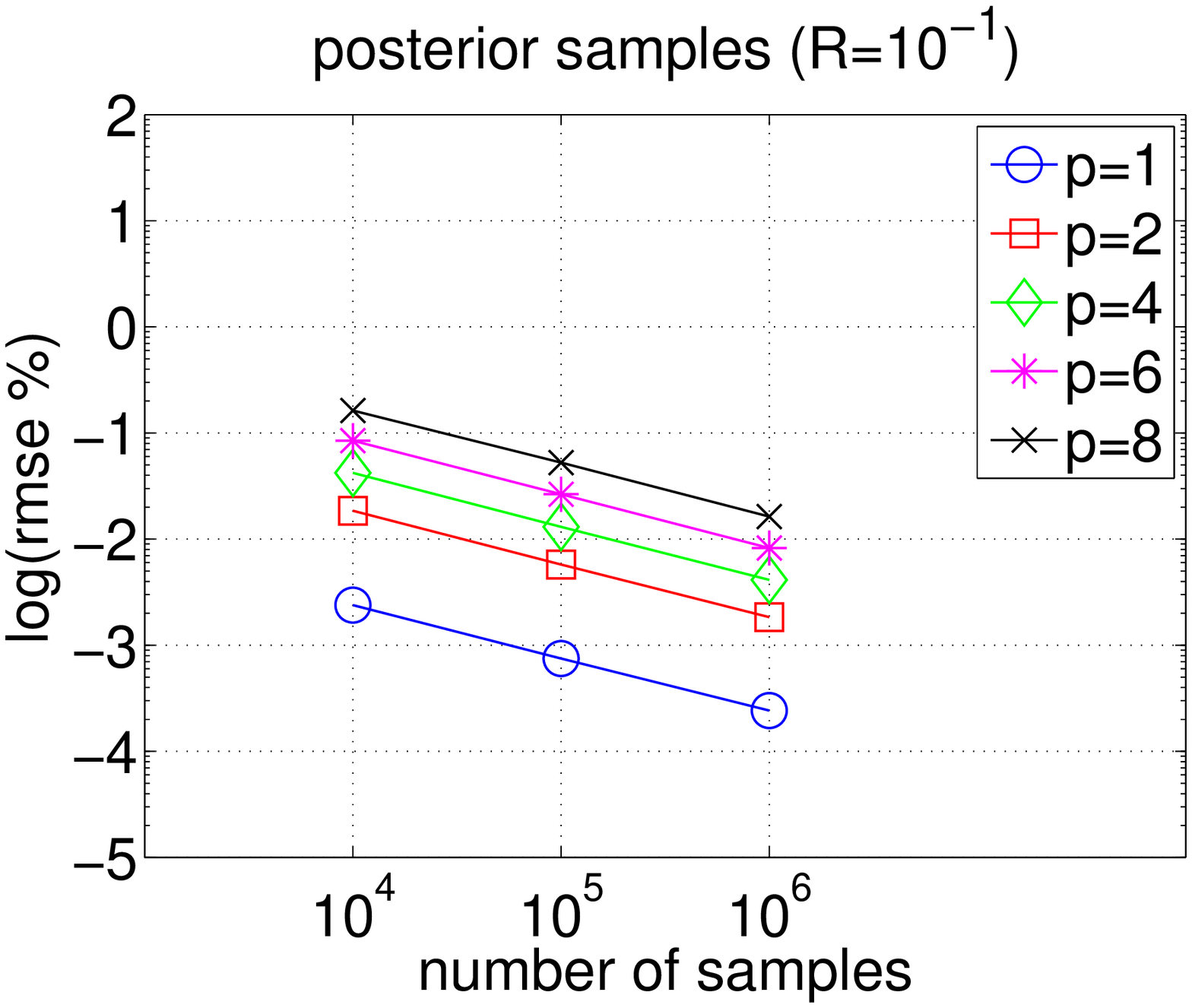} \label{fig:Da} } 
\subfigure[bootstrap reweighted prior samples]
{\includegraphics[width=0.39\textwidth]{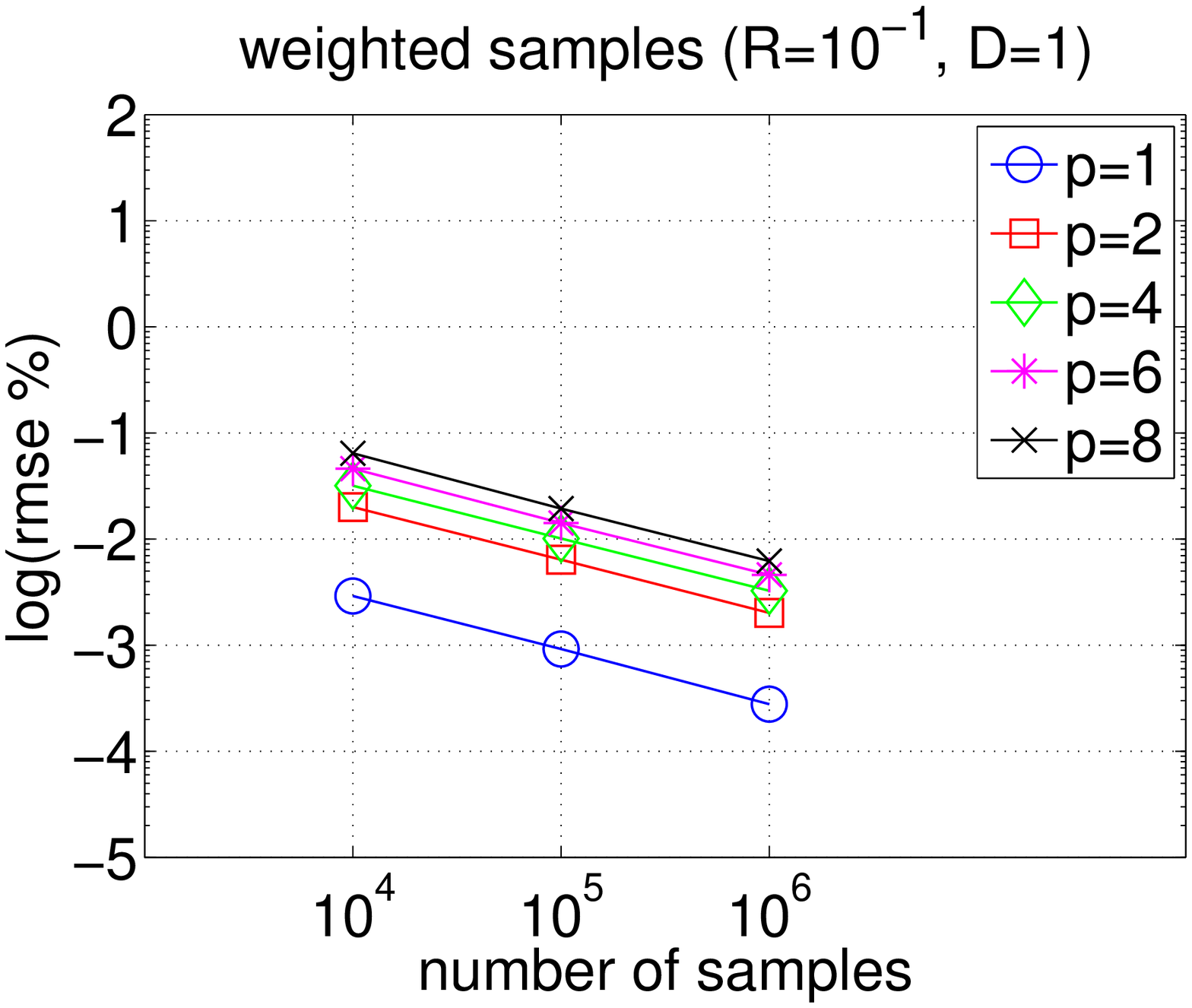} \label{fig:Db} } 
\subfigure[cubature approximation of posterior]
{\includegraphics[width=0.39\textwidth]{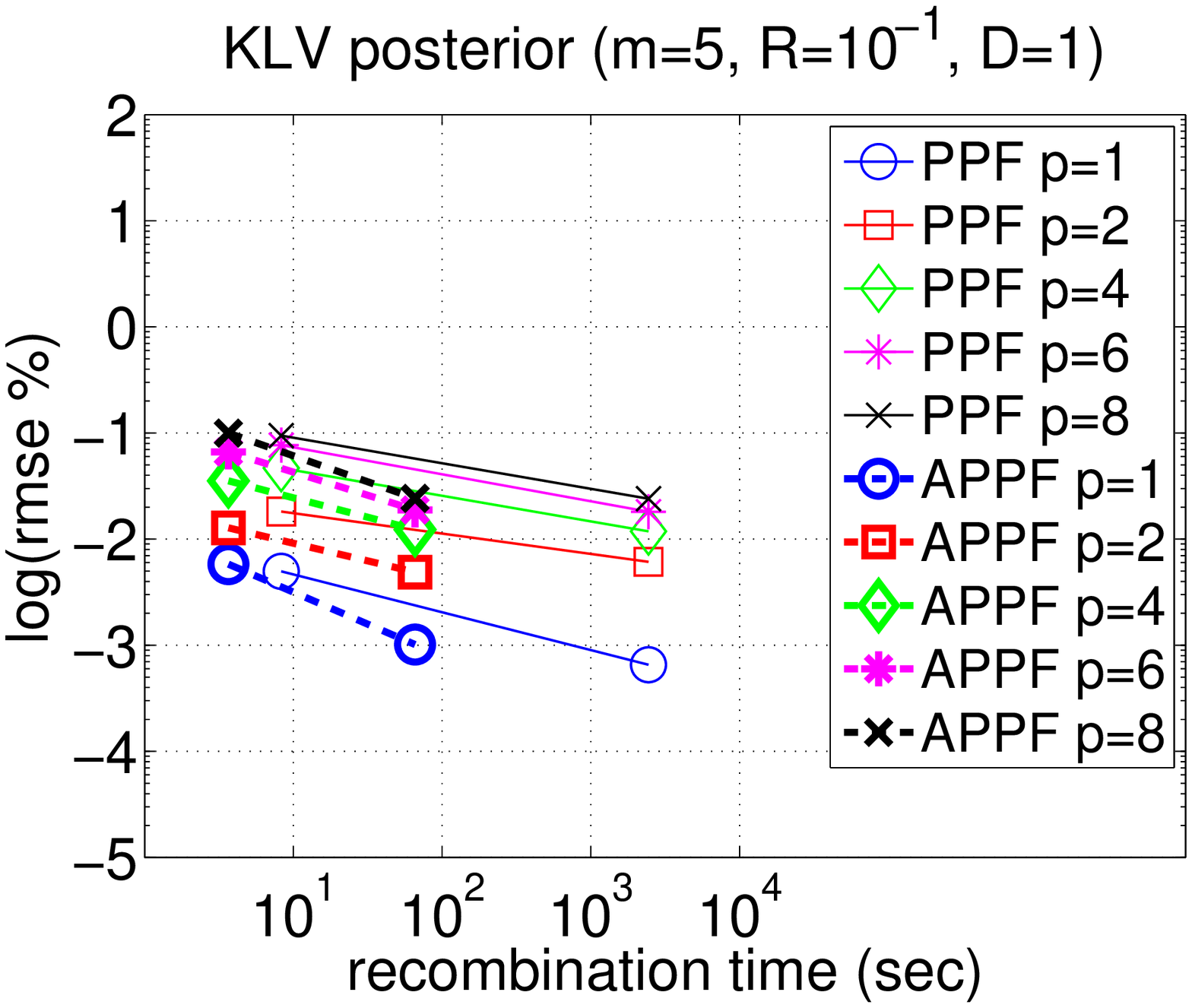} \label{fig:Dc}} 
} 
\centerline{
\subfigure[unweighted posterior samples]
{\includegraphics[width=0.39\textwidth]{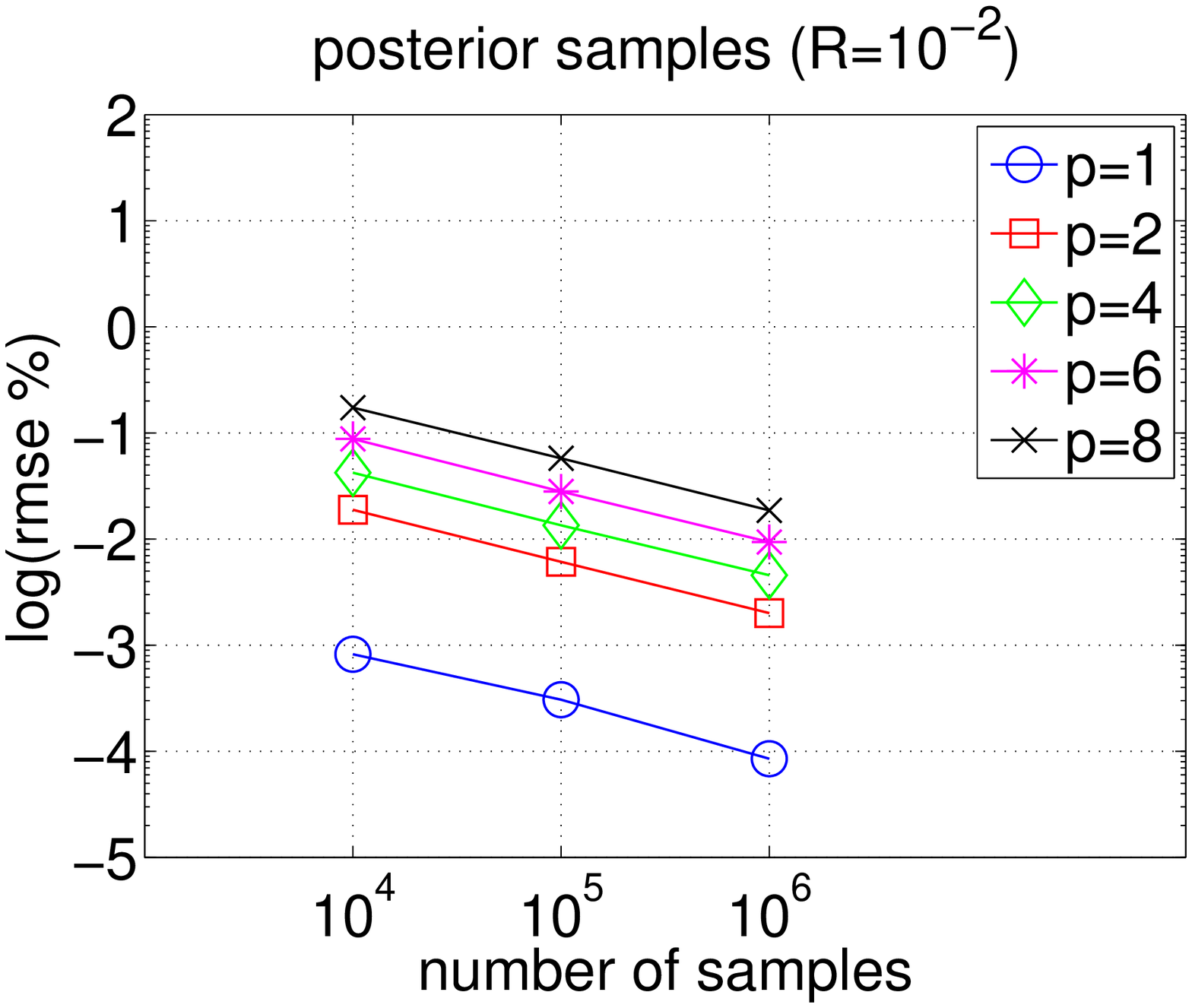} \label{fig:Dd}} 
\subfigure[bootstrap reweighted prior samples]
{\includegraphics[width=0.39\textwidth]{weightedpriorsamples_Rp01D1.eps} \label{fig:De}} 
\subfigure[cubature approximation of posterior]
{\includegraphics[width=0.39\textwidth]{cubpos_rp01d1.eps} \label{fig:Df}} 
}
\centerline{
\subfigure[unweighted posterior samples]
{\includegraphics[width=0.39\textwidth]{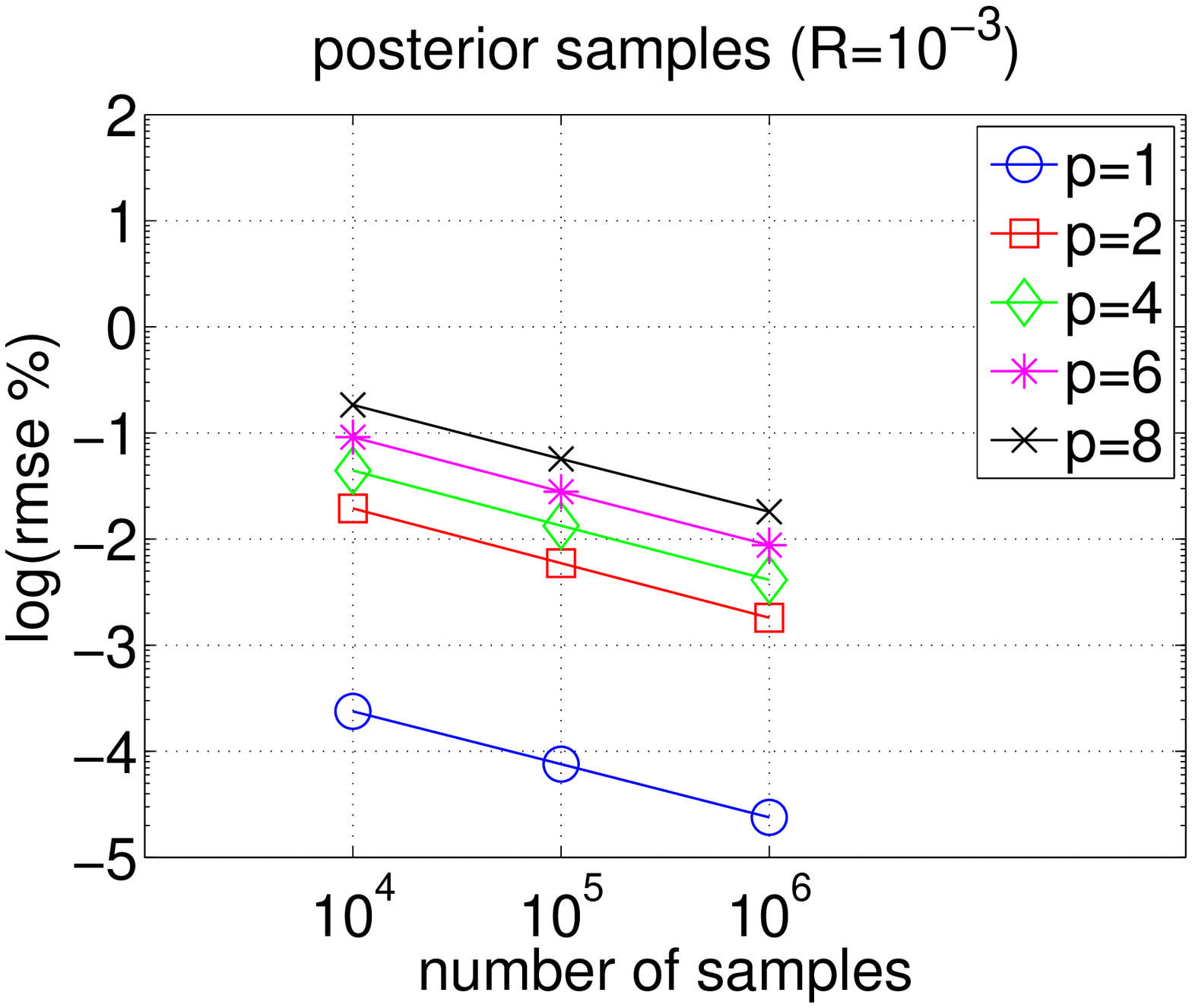} \label{fig:Dg}} 
\subfigure[bootstrap reweighted prior samples]
{\includegraphics[width=0.39\textwidth]{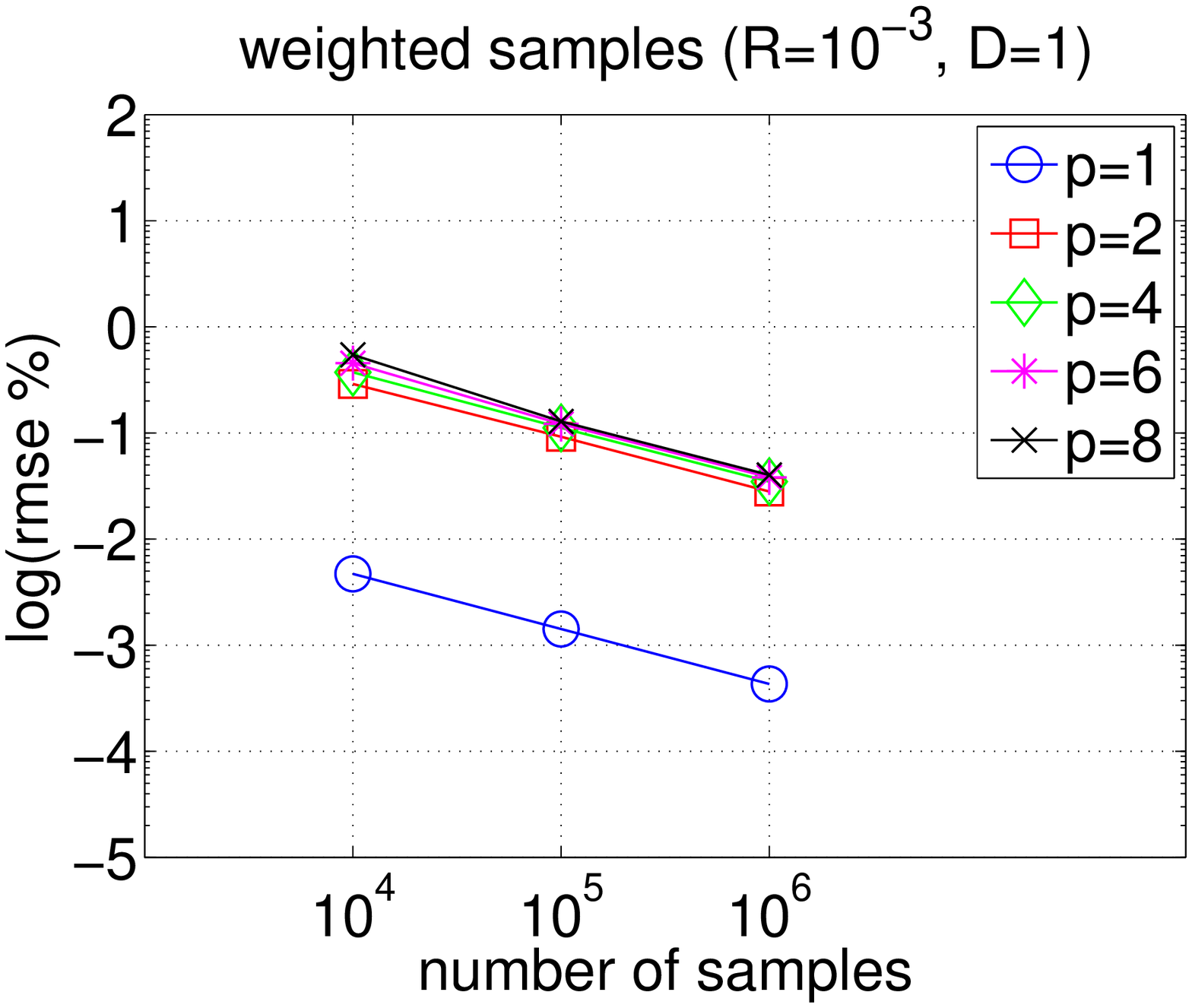} \label{fig:Dh}} 
\subfigure[cubature approximation of posterior]
{\includegraphics[width=0.39\textwidth]{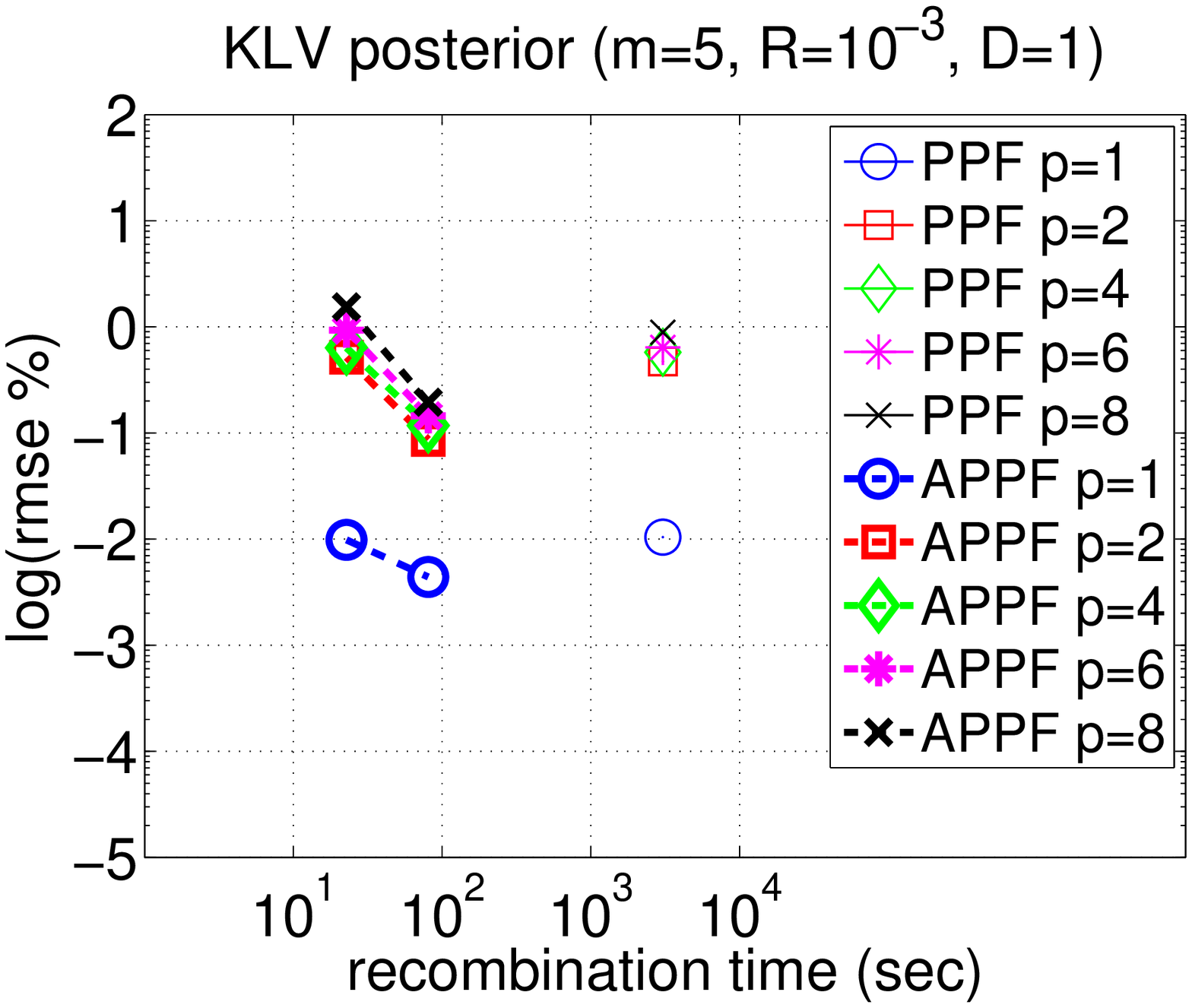} \label{fig:Di}} 
}
  \caption{
The posterior approximations when $D=1$ is fixed and $R=10^{-1},10^{-2},10^{-3}$ varies.
The left and middle columns are from Monte-Carlo samples
and the right column is form cubature approximation when $\epsilon=10^{-2}, 10^{-3}$.
In Fig.~\ref{fig:Di}, the PPF with $\epsilon = 10^{-3}$ is not produced.
} 
\label{fig:D} 
\end{figure}

In our numerical simulations,
the number of patches 
needed to satisfy
Eq.~(\ref{eq:supadarecomb})
in the PPF
increases as the time partition approaches the next observation time,
finally about $8^3 \sim 16^3$.
On the contrary,
Eq.~(\ref{eq:appfadarecomb}) 
in the APPF
is satisfied with
$2^3$ 
($< 10$)
patches in most cases.
As a result, APPF saves computation time significantly compared with PPF.

When $R=10^{-2}$ is fixed and $D=1,2,3$ varies,
the relative $L^2$ errors of the $p$-th moments of PPF and APPF
are shown in
Figs.~\ref{fig:Rb},
\ref{fig:Rd},
\ref{fig:Rf},
\ref{fig:Rh}.
We examine two cases of
$\epsilon=10^{-2}$
and 
$\epsilon=10^{-3}$.
The recombination times are measured
using Visual Studio
with Intel $2.53$ GHz processor
(the autonomous ODEs are solved analytically).
Fig.~\ref{fig:R} reveals the following.
\begin{itemize}
  \item
	The prior approximation of PPF with $\epsilon=10^{-3}$ shows similar
	accuracy with $10^4$ Monte-Carlo sampling (Figs.~\ref{fig:Ra}, \ref{fig:Rb}).
  \item The accuracy of the APPF prior approximation is in general worse than PPF especially for higher order
	moments (Fig.~\ref{fig:Rb}). 
  \item As the observation is located far from the prior mean,
i.e., as $D$ increases,
the posterior approximation obtained from Monte-Carlo bootstrap reweighting becomes 
less accurate 
(Figs.~\ref{fig:Rc}, \ref{fig:Re}, \ref{fig:Rg}).
\item The accuracy of the APPF posterior approximation
is similar to PPF but APPF significantly reduces the recombination time
(Figs.~\ref{fig:Rd}, \ref{fig:Rf}, \ref{fig:Rh}).
\item The accuracy of the PPF and APPF posterior approximations with $\epsilon=10^{-2}$
 is similar to $10^4$ Monte-Carlo reweighted samples when $D=1,2$ 
(Figs.~\ref{fig:Rc}, \ref{fig:Rd}, \ref{fig:Re}, \ref{fig:Rf})
 and to $10^5$ reweighted samples when $D=3$
(Figs.~\ref{fig:Rg}, \ref{fig:Rh}).
\item The accuracy of the PPF and APPF posterior approximations with $\epsilon=10^{-3}$
 is similar to $10^5$ Monte-Carlo reweighted samples when $D=1$
(Figs.~\ref{fig:Rc}, \ref{fig:Rd}),
to $10^6$ reweighted samples when $D=2$
(Figs.~\ref{fig:Re}, \ref{fig:Rf})
and 
to $10^7$ reweighted samples when $D=3$
(Figs.~\ref{fig:Rg}, \ref{fig:Rh}).
\end{itemize}

There is an important insight to be gained from this experimental analysis.
Though PPF
produces a better description of the
prior 
than APPF, 
this naive approximation of the prior 
is not good at
describing the 
fundamental object of interest in filtering, i.e., 
the posterior. 
The point is that one needs an
extremely accurate representation of the prior in certain localities. 
The 
APPF
delivers
this without undue cost. 
The 
PPF
method would have to deliver this accuracy uniformly and
well out into the tail of the prior.
As a result, PPF and APPF 
succeed in accurately describing the higher order statistics of posterior 
even when $D$ is big and Monte-Carlo fails to do such a job.

When $D=1$ is fixed and $R = 10^{-1},10^{-2},10^{-3}$ varies,
the values of 
Eq.~(\ref{eq:rmse})
for PPF and APPF
are shown in
Figs.~\ref{fig:Dc}, \ref{fig:Df}, \ref{fig:Di}.
We again examine two cases of
$\epsilon=10^{-2}$
and 
$\epsilon=10^{-3}$.
Fig.~\ref{fig:D} 
reveals the following.
\begin{itemize}
  \item 
	The moment approximations of Monte-Carlo Gaussian samples are insensitive to the covariance
	(recall that the diagonal element of $C_{n|n}$ is almost the same with the value of $R$)
	except the mean 
(Figs.~\ref{fig:Da}, \ref{fig:Dd}, \ref{fig:Dg}).
\item As the observation is more accurate,
i.e., as $R$ decreases,
the posterior approximation obtained from Monte-Carlo bootstrap reweighting becomes 
less accurate 
(Figs.~\ref{fig:Db}, \ref{fig:De}, \ref{fig:Dh}).
\item 
	As $R$ decreases,
  the recombination time 
  needed
  to achieve a given degree of accuracy
  becomes bigger for PPF
  but this is not the case for APPF, i.e., 
  the recombination time for APPF is insensitive to $R$
(Figs.~\ref{fig:Dc}, \ref{fig:Df}, \ref{fig:Di}).
\end{itemize}

The simulations show that the APPF becomes more competitive than the PPF
for the solution of the intermittent data assimilation problem
with small observation noise error.
It further shows that the APPF is of higher order with respect to the recombination time
and can achieve the given degree of accuracy with lower computational cost.

Although $Y_n$ is ``there and measurable'' it is sometimes the case 
that it is actually computationally very expensive to compute and 
that actually the thing one can compute is the evaluation of likelihood for a number of locations.
For example, consider a tracking problem for an object of moderate intensity 
and diameter that does a random walk and is moving against a slightly noisy background 
and is observed relatively infrequently. Its influence is entirely local. The likelihood function 
will be something like the Gaussian centred at the position of object but completely uninformative 
elsewhere in the space. The smaller the object, the tighter or narrower the Gaussian the harder 
the problem of finding the object becomes. 
One can compute the likelihood at any point in the space, 
but only evaluations at the location of the particle are informative. In that way one sees that
\begin{enumerate}
  \item 
The $Y_n$ is ``observable'' but only partially observed - and with low noise is very expensive to
observe accurately as one has to find the particle.
\item The likelihood can be observed at points in the space.
\end{enumerate}
In this sort of example it would be quite wrong to assume that, 
if we know the prior distribution of $X_n$ then  just because $Y_n = X_n + \eta_n$ 
we know the posterior distribution at zero cost. 
For sequential Monte-Carlo methods, 
bootstrap reweighting would seem to give a much better approach.

\subsection{Prospective performance PPF and APPF with cubature on Wiener space of degree $7$}
\label{sec:ppfappfd7}
A cubature formula on Wiener space of degree $m=7$ is currently not 
available 
when $d=3$.
Here we use another operator 
to see the prospective performance
of higher order cubature formula.

\begin{figure}
\centerline{
\subfigure[bootstrap reweighted prior samples]
{\includegraphics[width=0.39\textwidth]{weightedpriorsamples_Rp01D2.eps} \label{fig:D2a}} 
\subfigure[cubature approximation of posterior]
{\includegraphics[width=0.39\textwidth]{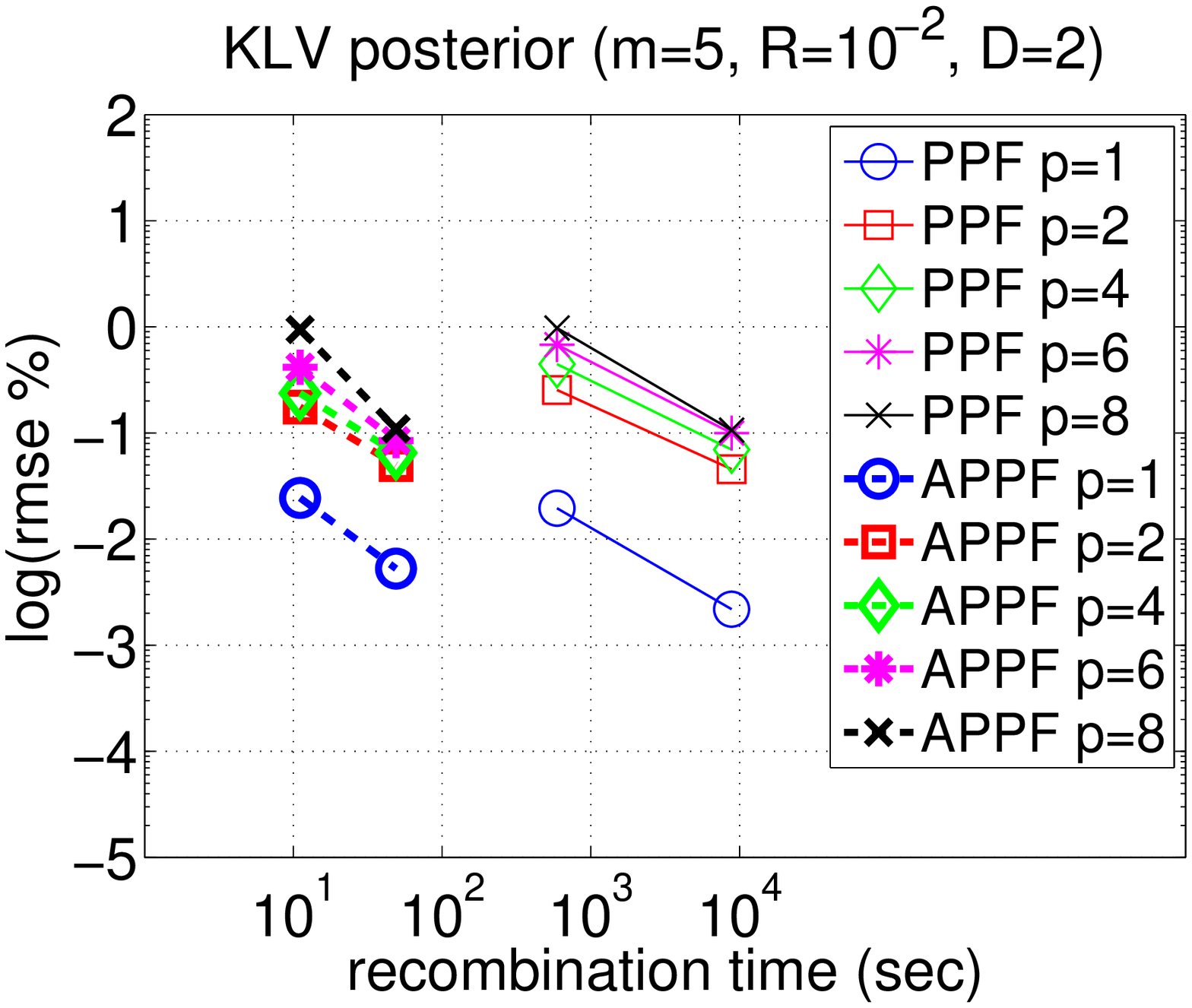} \label{fig:D2b}} 
\subfigure[cubature approximation of posterior]
{\includegraphics[width=0.39\textwidth]{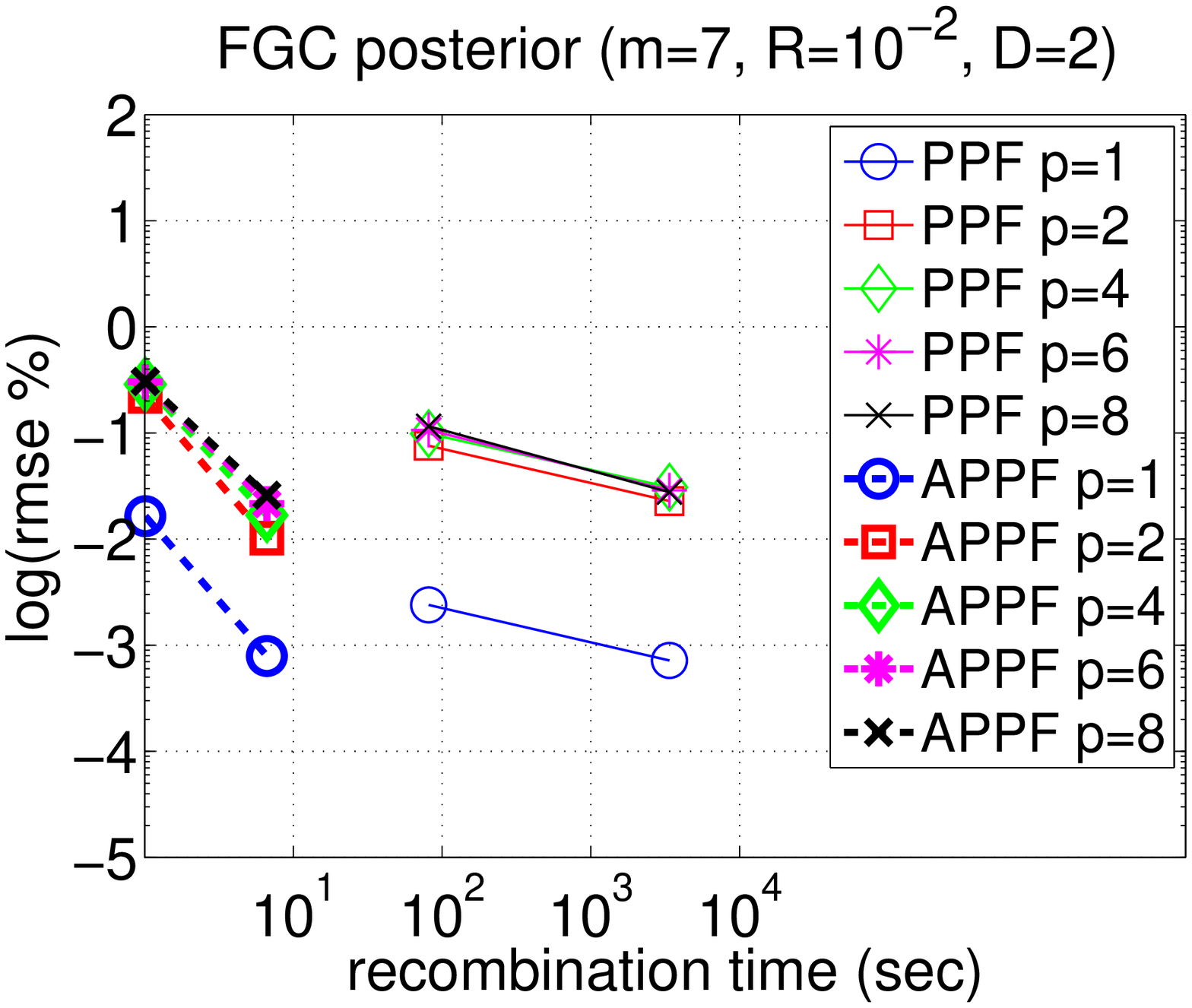} \label{fig:D2c}} 
} 
  \caption{
The posterior approximations when $R=10^{-2}$ and $D=2$.
The left is from Monte-Carlo samples 
and the middle and right is from cubature approximation when $\epsilon=10^{-2}, 10^{-3}$.
} 
\label{fig:D2} 
\end{figure}

\begin{table}
  \renewcommand{\arraystretch}{1.3} 
  \caption{The number of adaptive partition $k$ for FGC with $m=7$} 
  \label{tab:partm7} 
  \centering \begin{tabular}{ c |  c  c cc} 
	& $\epsilon = 10^{-2}$ & $\epsilon =10^{-3}$ & $\epsilon =10^{-4}$ & $\epsilon =10^{-5}$  \\ \hline
	$R=10^{-1}$  & 2 & 4 & 6 & 10  \\ 
	$R=10^{-2}$  & 5 & 9 & 16 & 28  \\ 
	$R=10^{-3}$  & 9 & 17 & 30 & 54  \\ 
  \end{tabular} 
\end{table}

For the linear dynamics
satisfying
\begin{equation*}
X(\Delta)=F_\Delta X(0) + \nu_\Delta, \quad \nu_\Delta \sim \mathcal{N}\left( 0, Q_\Delta \right)
\end{equation*}
for which $F_\Delta \in \mathbb{R}^{3\times 3}$ is a matrix,
one can define the forward operator 
\begin{equation}
  \label{eq:FGC}
{\text{FGC}}^{(m)}
\left(\sum_{i=1}^n \kappa_i \delta_{x^i}, \Delta \right)
\equiv \sum_{i=1}^n \sum_{j=1}^{n_m} \kappa_i \lambda_j \delta_{ F_\Delta x^i+ z^j}
\end{equation}
where $\{\lambda_j,z^j\}_{j=1}^{n_m}$ is a Gaussian cubature of degree $m$
with respect to the law of $\nu_\Delta$.
The authors see that
the performance of FGC 
is similar to KLV on the flow level when $m=3, 5$
and that
Eq.~(\ref{eq:FGC})
can be used as an alternative to 
Eq.~(\ref{eq:KLVflow})
for the PPF and APPF application to the test model, i.e.,
Eq.~(\ref{eq:L63}) where $a_0=0$.

The number of iterations $k$ in the adaptive partition,
obtained from using 
FGC
with
Gauss-Hermite cubature of degree $m=7$ whose support size is
$n_m = 64$
in place of $Q^{m}_{s_j}$,
is shown in 
table~\ref{tab:partm7}.
We apply FGC with degree $m=7$ 
to obtain a prior and posterior approximation,
where
the recombination degree $r=5$ and $\theta = 0.2\times \epsilon$ is used.
Our choice of $\tau$ is again such that $1/4 \sim 1/3$ of the particles
are allowed to skip to the next observation time.
Eq.~(\ref{eq:rmse})
in the case of $R=10^{-2}$, $D=2$ and $\epsilon=10^{-2}, 10^{-3}$
is shown
in Fig.~\ref{fig:D2c},
that can be viewed as 
an accuracy 
of PPF and APPF with cubature on Wiener space of degree $m=7$.
Its performance is in fact one higher order improvement for both accuracy and recombination time
in view of Figs.~\ref{fig:D2a}, \ref{fig:D2b}.
This result 
highlights the necessity to find cubature formula on Wiener space of degree $m=7$ 
in order to solve the PDE or filtering problem
with high accuracy in a moderate dimension.

\section{Discussion}
\label{sec:discussion} 
In this paper we introduce a hybrid
methodology for the numerical resolution of the filtering problem which we
call the adaptive patched particle filter (APPF). We explore some of its
properties and we report on a first attempt at a practical implementation.
The APPF combines many different \textquotedblleft methods\textquotedblright, 
each of which addresses a different part of the problem and has
independent interest. At a fundamental level all of the methods use high
order approaches to quantify uncertainty (cubature), and also to reduce the
complexity of calculations (recombination based on heavy numerical linear
algebra), while retaining explicit thresholds for accuracy in the individual
computation. The thresholds for accuracy in a stage are normally achievable
in a number of ways (e.g., small time step with low order, or large time step
with high order) and the determination of these choices depends on
computational cost. Aside from this use of the error threshold and choices
based on computational efficiency there are several other points to observe
in our development of this filter.

\begin{enumerate}
\item One feature is the surprising ease with which one can adapt the
computations to the observational data and so avoid performing unnecessary
computations. In even moderate dimensions (we work in $3+1$) this has a huge
impact for the computation time while preserving the accuracy we achieve for
the posterior distribution 
(Figs.~\ref{fig:Rd}, \ref{fig:Rf}, \ref{fig:Rh}, Figs.~\ref{fig:Dc}, \ref{fig:Df}, \ref{fig:Di}). 
It is an automated form of high order importance
sampling which has wider application than the one explored in this paper,
for instance
it is used to deliver
accurate answers to PDE problems with piecewise smooth 
test function
in the example developed in
\cite{litterer2012high}.
\item Another innovation allowing a huge reduction in computation is the
ability to efficiently \textquotedblleft patch the
particles\textquotedblright\ in the multiple dimensional scenario. Although
the problem might at first glance seem elementary, it is in fact the problem
of data classification. To resolve this problem we introduce an efficient
algorithm for data classification based on extending the Morton order to
floating point context. This method has now also been used effectively for
efficient function extrapolation 
\cite{wei}.
\item The KLV algorithm is at the heart of a number of successful methods
for solving PDEs in moderate dimension 
\cite{ninomiya2008weak}.
In each
case, something has to be done about the explosion of scenarios after each
time step; this in turn has to rely on and understanding the errors. In this
paper we take a somewhat different approach to the literature
\cite{lyons2004cubature}
in the way
we use higher order Lipschitz norms systematically to understand how well
functions have been smoothed, and to measure the scales on which they can be
well approximated by polynomials. This has the consequence that one can be
quite precise about the errors one incurs at each stage in the calculation.
In the end this is actually quite crucial to the logic of our approach since
an efficient method requires optimisation over several parameters -
something that is only meaningful if there are (at least in principle)
uniform estimates on errors. As a result of this perspective, we do not
follow the time steps and analytic estimates introduced 
by Kusuoka
in
\cite{kusuoka2004approximation}
although
we remain deeply influenced by balancing the smoothing properties of the
semigroup with the use of non-equidistant time steps.
\item The focus on Lipschitz norms makes it natural to apply an adaptive
approach to the recombination patches as well as to the prediction process.
In both cases we can be lead by the local smoothness of the likelihood
function as sampled on our high order high accuracy set of scenarios.
\item We have focussed our attention on the quality of the tail distribution
of the approximate posterior we construct. This is important in the
filtering problem because a failure to describe the tail behaviour of the
tracked object implies that one will lose the trajectory all together at
some point. These issues are particularly relevant in high dimensions as the
cost of increasing the frequency of observation can be prohibitive. If one
wishes to ensure reliability of the filter in the setting where there is a
significant discrepancy between the prior estimate and the realised outcome
over a time step then our APPF with cubature on Wiener space of degree $5$
already shows in the three dimensional example that it can completely
outperform sequential importance resampling Monte-Carlo approach. The
absence, at the current time, of higher order cubature formulae is in this
sense very frustrating as the evidence we give suggests that higher degree
methods will lead to substantial further benefits for both computation and
accuracy.
\end{enumerate}

In putting this paper together we have realised that there are many branches
in this algorithm that can be improved, in particular some parts of the
adaptive process and also the recombination (a theoretical improvement in
the order of recombination has recently been discovered
\cite{maria}).
There are also large parts that can clearly be parallelised.
We believe that there is ongoing scope for increasing the performance of the
APPF.

\section*{Acknowledgment}
The authors would like to thank the following institutions for their
financial support of this research. Wonjung Lee : NCEO project NERC and King
Abdullah University of Science and Technology (KAUST) Award No.
KUK-C1-013-04. Terry Lyons : NCEO project NERC, ERC grant number 291244 and
EPSRC grant number EP/H000100/1. The authors also thank the Oxford-Man
Institute of Quantitative Finance for its support.

\bibliographystyle{imsart-nameyear.bst}
\bibliography{bibliographyfile}

\end{document}